\documentclass[11pt]{amsart}
\usepackage{amsxtra}
\usepackage{amssymb}
\usepackage[mathscr]{eucal}

\theoremstyle{plain}
\newcommand{\cleqn}{\setcounter{equation}{0}}
\newcommand{\clth}{\setcounter{theorem}{0}}
\newcommand {\sectionnew}[1]{\section{#1}\cleqn\clth}

\theoremstyle{plain}
\newtheorem{theorem}{Theorem}[section]
\newtheorem{lemma}[theorem]{Lemma}
\newtheorem{definition-theorem}[theorem]{Definition-Theorem}
\newtheorem{proposition}[theorem]{Proposition}
\newtheorem{corollary}[theorem]{Corollary}

\theoremstyle{definition}
\newtheorem{definition}[theorem]{Definition}
\newtheorem{example}[theorem]{Example}

\newtheorem{remark}[theorem]{Remark}
\newtheorem{conjecture}[theorem]{Conjecture}
\newtheorem{notation}[theorem]{Notation}

\newtheorem{observations}[theorem]{Observations}
\newcommand \bth[1] { \begin{theorem}\label{t#1} }
\newcommand \ble[1] { \begin{lemma}\label{l#1} }

\newcommand \bpr[1] { \begin{proposition}\label{p#1} }
\newcommand \bco[1] { \begin{corollary}\label{c#1} }
\newcommand \bde[1] { \begin{definition}\label{d#1}\rm }
\newcommand \bex[1] { \begin{example}\label{e#1}\rm }
\newcommand \bre[1] { \begin{remark}\label{r#1}\rm }
\newcommand \bcj[1] { \begin{conjecture}\label{j#1}\rm }

\newcommand \bnota[1] { \begin{notation}\label{n#1}\rm }
\renewcommand {\eth} { \end{theorem} }
\newcommand {\ele} { \end{lemma} }

\newcommand {\epr} { \end{proposition} }
\newcommand {\eco} { \end{corollary} }
\newcommand {\ede} { \end{definition} }
\newcommand {\eex} { \end{example} }
\newcommand {\ere} { \end{remark} }
\newcommand {\ecj} { \end{conjecture} }

\newcommand {\enota} { \end{notation} }
\newcommand \thref[1]{Theorem \ref{t#1}}
\newcommand \leref[1]{Lemma \ref{l#1}}
\newcommand \prref[1]{Proposition \ref{p#1}}
\newcommand \coref[1]{Corollary \ref{c#1}}

\newcommand \deref[1]{Definition \ref{d#1}}
\newcommand \exref[1]{Example \ref{e#1}}

\def \Rset {{\mathbb R}}         
\def \Cset {{\mathbb C}}
\def \KK {{\mathbb K}}

\def \Zset {{\mathbb Z}}

\def \AA {{\mathcal{A}}}

\def \TT {{\mathcal{T}}}

\def \OO {{\mathcal{O}}}
\def \PP {{\mathcal{P}}}

\def \UU {{\mathcal{U}}}

\def \TT {{\mathcal{T}}} 

\def \bb {{\bf{b}}}
\def \cb {{\bf{c}}}
\def \db {{\bf{d}}}
\def \qb {{\bf{q}}}

\def \de {\delta}
\def \al {\alpha}
\def \be {\beta}

\def \la {\lambda}

\def \ga {\gamma}
\def \de {\delta}
\def \Ga {\Gamma}

\def \sig {\sigma}

\def \sig{\sigma}

\def \sy  {\ast}                         

\def \rcor {\rangle}
\def \lcor {\langle}

\def \ol {\overline}
\def \wt {\widetilde}

\def \id { {\mathrm{id}} }

\def \g  {\mathfrak{g}}   

\def \sl {\mathfrak{sl}}

\def \sl {\mathfrak{sl}}

\DeclareMathOperator \Aut { {\mathrm{Aut}} }

\DeclareMathOperator \gr  { {\mathrm{gr}} }
\DeclareMathOperator \fl { {\mathrm{fl}} }

\DeclareMathOperator \supp { {\mathrm{supp}} }

\DeclareMathOperator \Fract { {\mathrm{Fract}} }

\renewcommand \max { {\mathrm{max}} }
\newcommand\kx{\KK^*}

\newcommand\HH{{\mathcal{H}}}
\newcommand\xh{X(\HH)}
\DeclareMathOperator \Spec {Spec}

\newcommand \Znn {\Zset_{\ge 0}}
\newcommand \Zpos {\Zset_{> 0}}

\newcommand \Hmax {\HH_{\max}}

\def\lab{{\boldsymbol \lambda}}

\newcommand{\cbsh}{\cb^\#}

\newcommand{\bfp}{{\boldsymbol p}}
\newcommand{\bfq}{{\boldsymbol q}}
\newcommand{\Obfqkn}{\OO_\bfq(\KK^n)}
\DeclareMathOperator{\autgr}{Aut_{\gr}}
\DeclareMathOperator{\autfl}{Aut_{\fl}}
\DeclareMathOperator\AD{AD}
\DeclareMathOperator\GL{GL}
\DeclareMathOperator\SL{SL}
\DeclareMathOperator\tw{tw}
\DeclareMathOperator\img{im}
\newcommand\Ahat{\widehat{A}}
\newcommand\ghat{\widehat{g}}
\DeclareMathOperator \udim {u{.}dim}
\newcommand \kbar {\overline{\KK}}

\newcommand \OlpMn {\OO_{\la,\bfp}(M_n(\KK))}
\DeclareMathOperator\PIdeg{PI-deg}
\DeclareMathOperator\card{card}

\begin{document}

\title{Twist invariants of graded algebras}

\author[K. R. Goodearl]{K. R. Goodearl}
\address{
Department of Mathematics \\
University of California\\
Santa Barbara, CA 93106 \\
U.S.A.
}
\email{goodearl@math.ucsb.edu}
\author[M. T. Yakimov]{M. T. Yakimov}
\thanks{The research of K.R.G. was partially supported by NSF grant DMS-1601184,
and that of M.T.Y. by NSF grant DMS-1601862.}
\address{
Department of Mathematics \\
Louisiana State University \\
Baton Rouge, LA 70803 \\
U.S.A.
}
\email{yakimov@math.lsu.edu}
\date{}
\begin{abstract}
We define two invariants for (semiprime right Goldie) algebras, one for algebras graded by arbitrary abelian groups, which is unchanged under twists by  $2$-cocycles on the grading group, and one for $\Zset$-graded or $\Znn$-filtered algebras. The first invariant distinguishes quantum algebras which are 
``truly multiparameter" apart from ones that are ``essentially uniparameter", meaning  cocycle twists of uniparameter algebras. We prove that both invariants 
are stable under adjunction of polynomial variables. Methods for computing these invariants for large families of algebras are given, including 
quantum nilpotent algebras and algebras admitting one quantum cluster, and applications to non-isomorphism theorems are obtained.
\end{abstract}

\maketitle

\sectionnew{Introduction}

We introduce and investigate two interrelated invariants for graded or filtered algebras. The first concerns algebras $A$ graded by abelian groups $\Ga$, and takes account of the related algebras obtained by twisting $A$ with respect to $2$-cocycles on $\Ga$. The resulting invariant, which we denote $\tw_\Ga(A)$, is a cocycle twist invariant, and helps to separate algebras which are not twist equivalent. We also apply it to distinguish between quantum algebras which are basically uniparameter -- meaning the scalars involved in commutation relations reduce to powers of a single scalar -- and those which are truly multiparameter, in that their commutation scalars under any cocycle twist do not form a cyclic group. Our second invariant concerns $\Zset$-graded or $\Znn$-filtered algebras $A$ and applies the invariant $\tw_\Ga$ to a suitably maximal grading group $\Ga$, obtained as the character group of a maximal torus in an affine algebraic group of graded or filtered automorphisms. The resulting invariant, denoted $\tw(A)$, is used to help separate algebras which are not isomorphic as $\Zset$-graded or $\Znn$-filtered algebras, respectively.
These invariants, and a more basic invariant of Alev and Dumas \cite{AD} on which they are founded, are in many circumstances stable under polynomial extensions. Consequently, they allow to prove stability results saying that certain algebras cannot be made isomorphic even after adjoining arbitrarily many polynomial indeterminates.

Our original motivation for these invariants was to attack the question ``uniparameter or multiparameter?" for classes of quantum algebras, such as quantum nilpotent algebras.
We do so  by adapting an invariant that was originally introduced by Alev and Dumas in \cite{AD} to tell quantum affine spaces and their division rings of fractions apart. We put their invariant into the following context: For any semiprime right Goldie algebra $A$, we define $\AD(A)$ to be the intersection of $\kx$ with the commutator subgroup of the group of units of the right Goldie quotient ring $\Fract A$. If, additionally,  $A$ is graded by an abelian group $\Ga$, we define the ``$\Ga$-twist invariant" $\tw_\Ga(A)$ to be the intersection of the AD-invariants of all cocycle twists of $A$. 
For instance, consider a quantum affine space $A = \OO_\bfq(\KK^n)$, for some multiplicatively skew-symmetric matrix $\bfq = (q_{ij}) \subseteq M_n(\kx)$. In this case, $\AD(A)$ is the subgroup of $\kx$ generated by all $q_{ij}$ \cite[Proposition 3.9]{AD}. However,  when $A$ is given its natural $\Zset^n$-grading, $\tw_{\Zset^n}(A)$ is the trivial group, due to the fact that $A$ can be written as a cocycle twist of a commutative polynomial ring. This is an extreme case, since quantum affine spaces and their graded subalgebras are essentially the only quantum algebras that can be obtained from commutative algebras via cocycle twists. In the case of multiparameter quantum matrix algebras $A = \OlpMn$, the invariant $\AD(A)$ again records all the relevant commutation scalars ($\la$ and the entries of $\bfp$) whereas $\tw_\Ga(A)$ is the cyclic group generated by $\la$ (assuming $A$ is graded by $\Ga = \Zset^{2n}$ in the natural way); see \exref{OlpMn2}. On the other hand, the multiparameter quantized Weyl algebras $A = A_n^{Q,P}(\KK)$ are typically truly multiparameter: with respect to the natural $\Zset^n$-grading, $\tw_{\Zset^n}(A)$ is the group generated by the entries of the vector $Q$ (\exref{qWeyl}).

From the perspective of $\Zset$-graded algebras, the twist invariant may not provide a lot of new information, due to the paucity of ``new commutators" when twisting by a cocycle on $\Zset$ -- in particular, multiplicative commutators of homogeneous elements do not change under cocycle twists in this setting. However, there are many algebras for which $\Zset$-gradings can be extended to larger grading groups. We develop a procedure for obtaining gradings by abelian groups which are maximal in a suitable sense. Specifically, if an algebra $A$ is affine (finitely generated as an algebra) and has a locally finite dimensional $\Zset$-grading, the group of graded algebra automorphisms of $A$ has a natural structure of affine algebraic group, and $A$ receives a grading from the character group $\xh$ of any maximal torus $\HH$ of $\autgr A$. We define the ``general twist invariant" of $A$ to be $\tw(A) = \tw_{\xh}(A)$ which, as we observe, is independent of the choice of $\HH$ (\leref{tw.indepH}). A similar procedure allows to define $\tw(A)$ for algebras $A$ equipped with a locally finite dimensional $\Znn$-filtration, and this coincides with the previous definition if $A$ is $\Zset$-graded and the filtration is the natural one associated with the $\Zset$-grading (\leref{maxtor.autgr}). 

We prove that $\tw$ is unchanged in the passage from a $\Zset$-graded algebra $A$ to a polynomial ring $A[x]$, where the grading on $A$ is extended to $A[x]$ so that $x$ is homogeneous with degree $1$, and similarly for $\Znn$-filtered algebras (\thref{twA[x]}). This allows to show that certain algebras cannot become isomorphic after making polynomial extensions. For instance, consider quantized Weyl algebras $A = A_n^{Q,P}(\KK)$ and $A' = A_{n'}^{Q',P'}(\KK)$. If $\langle q_1, \dots, q_n \rangle \ne \langle q'_1, \dots, q'_n \rangle$, then no cocycle twisted 
polynomial algebra $\big( A[X_1, \dots, X_s] \big)^{\cb}$ is isomorphic, as a filtered $\KK$-algebra, to any cocycle twisted polynomial algebra 
$\big( A'[X_1, \dots, X_t] \big)^{\cb'}
$ (\exref{qWeyl3}). 
Here $A$ is filtered so that the standard generators $x_i$, $y_j$ for $A_n^{Q,P}(\KK)$ and the polynomial variables $X_k$ all have degree $1$.
\medskip 
\\
{\bf{Notation}}
Fix a base field $\KK$ throughout. Starting at the end of Subsection \ref{filtautgroups}, we will assume that $\KK$ is infinite, but otherwise $\KK$ is arbitrary.  We assume that all algebras are unital $\KK$-algebras, all automorphisms are $\KK$-algebra automorphisms, and that all skew derivations are $\KK$-linear. 

When we use the group $\Zset^n$, we denote its standard basis by $e_1,\dots,e_n$.
\sectionnew{$2$-cocycles and twists}
In this section we gather some background material on 2-cocycle twists of algebras. 
\bde{2co}
Let $\Ga$ be  an (additive) abelian group, and view $\kx$ as a $\Ga$-module with the trivial action. Then a \emph{$\kx$-valued $2$-cocycle on $\Ga$} is a map $\cb: \Ga \times \Ga \to \kx$ such that
\begin{equation}
\label{2co}
\cb(s, t+u) \cb(t,u) = \cb(s+t,u) \cb(s,t), \;\; \forall\, s,t,u \in \Ga.
\end{equation}
It is \emph{normalized} if $\cb(0,0) = 1$. Among these maps are the \emph{bicharacters} on $\Ga$, that is, the maps $\Ga\times\Ga \rightarrow \kx$ which are homomorphisms in each variable. A bicharacter $\cb$ is \emph{alternating} provided $\cb(s,s) = 1$ for all $s\in \Ga$, which implies skew-symmetry: $\cb(t,s) = \cb(s,t)^{-1}$ for all $s,t \in \Ga$.

Note that \eqref{2co} with $t = u = 0$ yields $\cb(s,0) \cb(0,0) = \cb(s,0)^2$. Thus, $\cb(s,0) = \cb(0,0)$ for all $s\in \Ga$, and similarly $\cb(0,s) = \cb(0,0)$. Consequently, if $\cb$ is normalized, then $\cb(s,0) = \cb(0,s) = 1$ for all $s\in \Ga$.

Denote by $Z^2(\Ga,\kx)$ the set of all $\kx$-valued $2$-cocycles on $\Ga$, and note that $Z^2(\Ga,\kx)$ is an abelian group under pointwise multiplication.
\ede

Recall that a matrix $(\al_{ij})$ of scalars is \emph{multiplicatively skew-symmetric} provided $\alpha_{ii} = 1$ and $\al_{ij}\alpha_{ji} = 1$ for all indices $i$, $j$.
 
\bpr{ASTprop}
Suppose $\Ga$ is free abelian with a basis $\{ e_i \mid i\in I \}$.

{\rm(a)} The rule $\rho(\cb) = \bigl( \cb(e_i, e_j) \cb(e_j, e_i)^{-1} \bigr)$ defines a surjection 
$$
\rho : Z^2(\Ga,\kx) \longrightarrow  
 \{ \text{multiplicatively skew-symmetric} \;\; I\times I \;\; \text{matrices over} \;\; \kx \}.
$$

{\rm (b)} For $\cb,\cb' \in Z^2(\Ga,\kx)$, one has $\rho(\cb) = \rho(\cb')$ if and only if there is a function $f: \Ga \rightarrow \kx$ such that
\begin{equation}
\label{c'=fc}
\cb'(s,t) = \frac{f(s) f(t)} {f(s+t)} \, \cb(s,t), \; \; \forall\, s,t \in \Ga.
\end{equation}

{\rm (c)} If $\cb \in Z^2(\Ga,\kx)$ and $\cbsh: \Ga \times \Ga \to \kx$ is defined by
$$
\cbsh(s,t) = \cb(s,t) \cb(t,s)^{-1}, \; \; \forall\, s,t \in \Ga,
$$
then $\cbsh$ is an alternating bicharacter on $\Ga$. 
\epr

\begin{proof} \cite[Proposition 1]{AST}.
\end{proof}

\bco{2cocycleform}
Suppose $\Ga$ is free abelian with a basis $\{ e_i \mid i\in I \}$, and put a total order on $I$. Let $\cb \in Z^2(\Ga,\kx)$, and assume there exists $p_{ij} = \sqrt{ \cb(e_i, e_j) \cb(e_j, e_i)^{-1} } \in \kx$ for all $i < j$. Let $\bb : \Ga \times \Ga \to \kx$ be the alternating bicharacter such that $\bb(e_i, e_j) = p_{ij}$ for all $i < j$. Then there is a function $f: \Ga \rightarrow \kx$ such that
$$
\cb(s,t) = \frac{f(s) f(t)} {f(s+t)} \, \bb(s,t), \; \; \forall\, s,t \in \Ga.
$$
\eco

\begin{proof} By construction, $\bb(e_i, e_j) \bb(e_j, e_i)^{-1} = p_{ij}^2 = \cb(e_i, e_j) \cb(e_j, e_i)^{-1}$ for all $i < j$. Since $\bb$ is alternating, it follows that $\rho(\bb) = \rho(\cb)$ for $\rho$ as in \prref{ASTprop}(a). Apply \prref{ASTprop}(b) to $\cb$ and $\bb$.
\end{proof}

\bde{twists}
Suppose $A$ is a $\Ga$-graded $\KK$-algebra and $\cb \in Z^2(\Ga,\kx)$. Define a bilinear multiplication operation $*_\cb$ on $A$ by the rule
\begin{equation}
\label{twistmult}
x *_\cb y := \cb(s,t) xy, \; \; \forall\, x \in A_s, \; \; y \in A_t, \; \; s,t \in \Ga.
\end{equation}
Then $(A,*_\cb)$ is a $\Ga$-graded $\KK$-algebra, called the \emph{twist of $A$ by $\cb$} \cite{AST}. We shall denote this algebra by $A^\cb$.
\ede

\begin{observations}  \label{twistobs}
Let $A$ be a $\Ga$-graded $\KK$-algebra.

(1) If $\cb \in Z^2(\Ga,\kx)$, then the identity element in $A^\cb$ is $\cb(0,0)^{-1} \cdot 1$, so if $\cb$ is normalized, this identity equals the identity of $A$. Moreover, in this case the subalgebras $A_0 \subseteq A$ and $A_0^\cb \subseteq A^\cb$ are identical.

(2) Suppose $\cb, \cb' \in Z^2(\Ga,\kx)$  and there is a function $f: \Ga \rightarrow \kx$ such that \eqref{c'=fc} holds.
Then $A^\cb$ and $A^{\cb'}$ are isomorphic, via the $\Ga$-graded  linear isomorphism given by 
$$
x \mapsto f(s)^{-1} x, \; \; \forall\, x \in A_s, \; \; s \in \Ga
$$
\cite{AST}.
Thus, starting with $\cb$, choosing $f$ such that $f(0) = \cb(0,0)^{-1}$, and defining $\cb'$ by \eqref{c'=fc}, we obtain $\cb'(0,0) = 1$. This yields a twist $A^{\cb'}$, isomorphic to $A^\cb$, in which the identity is the identity of $A$. 

(3) For any $\cb, \db \in Z^2(\Ga,\kx)$,
\begin{equation}  \label{doubletwist}
A^{\cb\db} = (A^\cb)^\db.
\end{equation}
Hence, the relation of \emph{cocycle twist equivalence}, where $\Ga$-graded $\KK$-algebras $A$ and $B$ are cocycle twist equivalent if and only if $B \cong A^\cb$ (as $\Ga$-graded algebras) for some $\cb \in Z^2(\Ga,\kx)$, is an equivalence relation.
\end{observations}

\sectionnew{Twist invariants for $\Gamma$-graded algebras}
In this section we review the Alev-Dumas invariant of noncommutative algebras and define a twist invariant of graded algebras. 
We also provide various examples to illustrate both notions.

\bde{defAD} Suppose $A$ is a semiprime right Goldie $\KK$-algebra. The \emph{Alev-Dumas invariant} of $A$ \cite[D\'efinition 3.8]{AD} is the group
$$\AD(A) := \GL_1(\Fract A)' \cap \kx,$$
where $\GL_1(\Fract A)'$ is the derived subgroup (= commutator subgroup) of the group of units of the classical right quotient ring $\Fract A$. 
\ede

\begin{observations}  \label{ADobs}
(1) If $A$ and $B$ are semiprime right Goldie $\KK$-algebras with $\Fract A \cong \Fract B$, then $\AD(A) = \AD(B)$.

(2) If $A_1,\dots,A_k$ are semiprime right Goldie $\KK$-algebras, then
\begin{equation}  \label{ADprod}
\AD(A_1\times \cdots\times A_k) = \AD(A_1) \cap \cdots\cap \AD(A_k).
\end{equation}

(3) For all positive integers $n$,
\begin{equation}  \label{ADMnK}
\AD(M_n(\KK)) = \{ \alpha \in \KK \mid \alpha^n=1 \}.
\end{equation} 
For $n\ge 3$, this follows from the well known fact that $\GL_n(\KK)' = \SL_n(\KK)$ (e.g., \cite[Theorem XIII.9.2]{Lang}). For $n=2$, it follows from that fact that $\GL_2(\KK)' \subseteq \SL_2(\KK)$ together with the observation 
$\left[\begin{smallmatrix} -1&0\\ 0&1\end{smallmatrix}\right] \left[\begin{smallmatrix} 0&1\\ 1&0\end{smallmatrix}\right] \left[\begin{smallmatrix} -1&0\\ 0&1\end{smallmatrix}\right] \left[\begin{smallmatrix} 0&1\\ 1&0\end{smallmatrix}\right] = \left[\begin{smallmatrix} -1&0\\ 0&-1\end{smallmatrix}\right]$.

Consequently, the $n$-th roots of unity in $K$ are contained in $\AD(M_n(A))$ for any semiprime right Goldie $\KK$-algebra $A$.

(4) The AD-invariant obviously depends on the choice of base field. For instance, if we view $M_3(\Cset)$ as a $\Cset$-algebra, then $\AD(M_3(\Cset))$ has order $3$, while as an $\Rset$-algebra, $\AD(M_3(\Cset))$ is the trivial group.
\end{observations}

For any multiplicatively skew-symmetric matrix $\bfq = (q_{ij}) \in M_n(\kx)$, the corresponding \emph{quantum affine space} is the $\KK$-algebra 
\[
\Obfqkn:= \KK \lcor y_1, \ldots, y_n \mid y_i y_j = q_{ij} y_j y_i \; \forall\, i,j \in [1,n] \rangle.
\]
The standard (one-parameter) quantum affine space, denoted $\OO_q(\KK^n)$ for $q\in \kx$, is the case of $\Obfqkn$ where $q_{ij} = q$ for $1\le i < j \le n$.
The localization of $\Obfqkn$ at the multiplicative subset generated by $y_1, \ldots, y_n$ is the \emph{quantum torus} 
\[
\TT^n_\bfq(\KK) = \OO_\bfq((\kx)^n) := \KK \lcor y_1^{\pm 1}, \ldots, y_n^{\pm 1} \mid y_i y_j = q_{ij} y_j y_i \; \forall\, i,j \in [1,n] \rangle.
\]
By way of abbreviation, we write $\langle \bfq \rangle$ for the group $\langle q_{ij} \mid i,j \in [1,n] \rangle \subseteq \kx$.

\bpr{ADqaff} {\rm\cite[Proposition 3.9]{AD}} Let $n \in \Zpos$ and $\bfq \in M_n(\kx)$ a multiplicatively skew-symmetric matrix. Then $\AD(\Obfqkn) = \langle \bfq \rangle$.
\epr

Of course, $\AD(\TT^n_\bfq(\KK)) = \langle \bfq \rangle$ as well. 

\bpr{PI} Let $A$ be a prime affine PI $\KK$-algebra of PI degree $n$.Then 
\[
\AD(A) \subseteq \{ \al \in \KK \mid \al^{n} =1 \}.
\]
\epr

\begin{proof} Denote by $Q$ the field of fractions of $Z := Z(A)$.  By Posner's theorem, $B : = A \otimes_Z Q$ is a central simple algebra over $Q$, and by definition, $n = \sqrt{\dim_Q B}$. This algebra splits over the algebraic closure $\ol{Q}$,
\[
A \otimes_Z \ol{Q} \cong M_n(\ol{Q}), 
\]
see \cite[Corollary 2.3.25]{Ro}. Therefore $A$ is isomorphic to a $\KK$-subalgebra of $M_n(\ol{Q})$, and by Observation \ref{ADobs}(1), 
\[
\AD(A) \subseteq  \AD(M_n(\ol{Q})) \cap \kx = \{ \al \in \KK \mid \al^{n} =1 \}.  \qedhere
\]
\end{proof}

Further examples appear in the following section.

\bde{deftwGam} Suppose $A$ is a $\Gamma$-graded semiprime right Goldie $\KK$-algebra. The \emph{twisted commutation invariant} of $A$ as a $\Gamma$-graded $\KK$-algebra is the group
$$\tw_\Gamma(A) := \bigcap_{\cb \in Z^2(\Gamma,\kx)} \AD(A^\cb) \subseteq \kx.$$
\ede

\bex{twGam.qaff} For any multiplicatively skew-symmetric matrix $\bfq = (q_{ij}) \in M_n(\kx)$, the quantum affine space $A = \Obfqkn$ has a natural grading by the group $\Ga = \Zset^n$ with $\deg y_i = e_i$ for all $i$.
This graded algebra can be
expressed as a cocycle twist
$$
A = \OO(\KK^n)^\cb,
$$
where $\OO(\KK^n) = \KK[x_1,\dots,x_n]$ is a polynomial algebra and $\cb$ is the $2$-cocycle such that
$$
\cb(s,t) = \prod_{1\le i<j\le n} q_{ij}^{s_it_j}, \;\; \forall\; s = (s_1,\dots,s_n),\, t = (t_1,\dots,t_n) \in \Ga.
$$
Since $\AD(\OO(\KK^n))$ is trivial (because $\OO(\KK^n)$ is commutative), we conclude that
$$
\tw_\Gamma(\Obfqkn) = \{ 1 \}.
$$
\eex

\bex{OlpMn} The \emph{$n\times n$ multiparameter quantum matrix algebra} $\OlpMn$ is given by generators $X_{ij}$ for $i,j\in [1,n]$ and relations
\begin{equation}  \label{OlpMn.rels}
X_{l m}X_{ij} = \begin{cases} 
p_{l i}p_{jm}X_{ij}X_{l m} +
(\lambda -1)p_{l i}X_{im}X_{l j}\,, &(l >i, \; m>j),  \\ 
\lambda p_{l i}p_{jm}X_{ij}X_{l m}\,, &(l >i, \; m\le j),  \\ 
p_{jm}X_{ij}X_{l m}\,, &(l=i, \; m>j),
\end{cases}
\end{equation}
where $\la \in \kx$ and $\bfp = (p_{ij}) \in M_n(\kx)$ is multiplicatively skew-symmetric. The \emph{standard $n\times n$ quantum matrix algebra} $\OO_q(M_n(\KK))$, for $q\in \kx$, is the case of $\OlpMn$ for which $\la = q^{-2}$ and $p_{ij} = q$ for all $i>j$. 

It is obvious from \eqref{OlpMn.rels} that $\AD(\OlpMn)$ contains $\la$ and all $p_{ij}$. We show below \eqref{ADtw.OlpMn} that if $\la$ is not a root of unity, then $\AD(\OlpMn) = \langle \la, \bfp \rangle$. In particular, $\AD(\OO_q(M_n(\KK))) = \langle q \rangle$ if $q$ is not a root of unity.

There is a standard grading of $\OlpMn$ by $\Ga = \Zset^n \times \Zset^n$, under which each $X_{ij}$ is homogeneous of degree $(e_i,e_j)$. If $\cb$ is the $2$-cocycle on $\Ga$ given by
$$
\cb\bigl( (s,t), (s',t') \bigr) = \prod_{1\le i < j \le n} p_{ij}^{s_j s'_i - t_j t'_i},
$$
then we find that
$$
\OlpMn^\cb \cong \OO_{\la,\mathbf1}(M_n(\KK)),
$$
where $\mathbf1$ denotes the matrix in $M_n(\kx)$ with all entries equal to $1$. Thus, by the result mentioned above, $\AD(\OlpMn^\cb) = \langle \la \rangle$, assuming that $\la$ is a non-root of unity. In fact, $\tw_\Ga(\OlpMn) = \langle \la \rangle$, as we show in \eqref{ADtw.OlpMn}.
\eex

\sectionnew{The twist invariant for algebras with one quantum cluster}
In this section we prove a general theorem that can be used to effectively compute the twist invariants of algebras that have one quantum cluster,
that is, algebras squeezed between a quantum affine space algebra and the corresponding quantum torus. The theorem is applied to compute the 
twist invariants of the multiparameter quantized Weyl algebras in full generality.

\subsection{The AD invariant for quantum cluster algebras}
A {\em{quantum cluster algebra}} $\AA$ is generated by a collection of cluster variables, whose different clusters are obtained from each other by mutation.
Each quantum seed of $\AA$ contains elements $y_1, \ldots, y_n  \in \AA$ such that $y_i y_j = q_{ij} y_j y_i$ for some $q_{ij} \in \kx$ and, in addition, 
these elements give rise to embeddings
\begin{equation}
\label{aa-emb}
\Obfqkn \subseteq \AA \subset \TT^n_\bfq(\KK) \,,
\end{equation}
where $\bfq = (q_{ij})$ is a multiplicatively skew-symmetric matrix.
We refer the reader to \cite{BeZ} for details on the uniparameter case when the scalars $q_{ij}$ are powers of a single scalar $q \in \kx$ and to
\cite[Ch. 2]{GoYa2} for the general case. To each quantum cluster algebra we associate the abelian group
\[
q(\AA):= \lcor  \bfq \rcor.
\]
This is independent of the choice of seed due to the mutation formula \cite[Eq. (2.19)]{GoYa2} for the matrices $\bfq$ of the quantum seeds of $\AA$. 
The {\em{upper quantum cluster algebra}} $\UU$ associated to $\AA$ (with or without inverted frozen variables, see \cite{BeZ,GoYa2}) also satisfies
\begin{equation}
\label{uu-emb}
\Obfqkn \subseteq \UU \subset \TT_\bfq^n(\KK).
\end{equation}
Set also $q(\UU):= q(\AA) =  \lcor  \bfq \rcor.$ Analogously to the case for $\AA$, this is independent of the choice of seed. 

A $\KK$-algebra $\AA$ that satisfies (possibly only one) embedding of the form \eqref{aa-emb} is called {\em{an algebra with a quantum cluster}}.
The class of algebras with one quantum cluster is strictly larger than that of quantum cluster algebras.

\bpr{AD-q-clu} {\em{(a)}} If $A$ is a $\KK$-algebra with a quantum cluster 
\[
\Obfqkn \subseteq A \subseteq \TT^n_\bfq(\KK) 
\]
for some multiplicatively skew-symmetric matrix $\bfq \in M_n(\kx)$, then $\AD(A) = \lcor \qb \rcor$.

{\em{(b)}} For every quantum cluster algebra $\AA$ and the associated upper quantum cluster algebra $\UU$,
\[
\AD(\AA) = \AD(\UU)= q(\AA).
\]
\epr
The statement follows from \prref{ADqaff}, Observation \ref{ADobs}(1) and the embeddings \eqref{aa-emb}-\eqref{uu-emb}.

Propositions \ref{pPI} and \ref{pAD-q-clu} have the following interesting  corollary. 

\bco{PI-qcl} Assume that $A$ is a $\KK$-algebra with a quantum cluster
\[
\Obfqkn \subseteq A \subseteq \TT^n_\bfq(\KK)
\]
for some multiplicatively skew-symmetric matrix $\bfq \in M_n(\kx)$.
Then $A$ is a PI algebra if and only if $\AD(A)$ is finite, in which case
$$
\PIdeg A \ge \card \AD(A)\,.
$$
 \eco 
 
 \begin{proof}
 It is well known that $\Obfqkn$ is a PI algebra if and only if all $q_{ij}$ are roots of unity (e.g., \cite[Theorem 7]{LM}), i.e., if and only if $\langle \bfq \rangle$ is finite. Since $\Obfqkn \subseteq A \subseteq \Fract \Obfqkn$, one of these algebras is PI if and only if the other is. The first result thus follows from \prref{AD-q-clu}(a).
 
If $A$ is PI with $\PIdeg A = d$, then $\card \AD(A) \le d$ by \prref{PI}.
\end{proof}

The inequality in \coref{PI-qcl} is by no means sharp, as can be seen from the formula for $\PIdeg \Obfqkn$ in \cite[Proposition 7.1(c)]{DeCP}. For instance, if $\bfq$ is given by $q_{ij} = q^{a_{ij}}$ where $q \in \kx$ is a primitive $\ell$-th root of unity and $(a_{ij}) \in M_n(\Zset)$ is an invertible antisymmetric matrix, then $\PIdeg \Obfqkn = \ell^{n/2}$.

\subsection{The twist invariant for algebras with one quantum cluster} \bde{Ochi}
Let $\chi : \Zset^n \times \Zset^n \rightarrow \kx$ be an alternating bicharacter. Set
$$
\bfq(\chi) := (\chi(e_i,e_j)) \in M_n(\kx),
$$
which is a multiplicatively skew-symmetric matrix, and
$$
\OO_\chi(\KK^n) := \OO_{\bfq(\chi)}(\KK^n), \qquad \TT^n_\chi(\KK) := \TT^n_{\bfq(\chi)}(\KK).
$$
\ede

The next theorem provides an effective way for computing twist invariants of algebras with one quantum cluster.

\bth{twGa.sandwich}
Let $\chi : \Zset^n \times \Zset^n \rightarrow \kx$ be an alternating bicharacter and
$$
\phi: \Zset^n \rightarrow \Ga
$$
a surjective group homomorphism, where $\Ga$ is a free abelian group. Assign to $\OO_\chi(\KK^n)$ and $\TT_\chi(\KK)$ the $\Ga$-gradings induced from their $\Zset^n$-gradings along $\phi$. For every inclusion of $\Ga$-graded $\KK$-algebras
$$
\OO_\chi(\KK^n) \subseteq A \subseteq \TT^n_\chi(\KK),
$$
we have
\begin{equation}  \label{twGa.eqn}
\tw_\Ga(A) = \langle \chi(\ker \phi, \, \Zset^n) \rangle \subseteq \kx.
\end{equation}
\eth

\begin{proof} For $\cb \in Z^2(\Ga,\kx)$, let $\cbsh$ be the alternating bicharacter on $\Ga$ obtained from $\cb$ as in \prref{ASTprop}(c), and define the alternating bicharacter $\chi_\cb$ on $\Zset^n$ by
\begin{equation}  \label{chi_c}
\chi_\cb(s,t) = \chi(s,t)\, \cbsh(\phi(s), \phi(t)).
\end{equation}
It is easy to see that there are compatible $\KK$-algebra isomorphisms
$$
\OO_\chi(\KK^n)^\cb \rightarrow \OO_{\chi_\cb}(\KK^n) \qquad\text{and}\qquad \TT^n_\chi(\KK)^\cb \rightarrow \TT^n_{\chi_\cb}(\KK).
$$
The second isomorphism sends $A^\cb$ to a subalgebra of $\TT^n_{\chi_\cb}(\KK)$ containing $\OO_{\chi_\cb}(\KK^n)$, and so \prref{AD-q-clu}(a) implies that $\AD(A^\cb) = \langle \bfq(\chi_\cb) \rangle$. Therefore
$$
\tw_\Ga(A) = \bigcap_{\cb \in Z^2(\Ga,\kx)} \langle \bfq(\chi_\cb) \rangle.
$$
Given $s\in \ker \phi$ and $t\in \Zset^n$, we see that $\chi_\cb(s,t) = \chi(s,t)$ for any $\cb \in Z^2(\Ga,\kx)$. Thus,
$$
\tw_\Ga(A) \supseteq \langle \chi(\ker \phi, \, \Zset^n) \rangle.
$$
We are left with proving the opposite inclusion.

Since $\Ga$ is free abelian, there is a basis $b_1,\dots,b_n$ for $\Zset^n$ such that some initial sublist $b_1,\dots,b_m$ is a basis for $\ker \phi$ and $\phi(b_{m+1}), \dots, \phi(b_n)$ is a basis for $\Ga$. Define $\cb \in Z^2(\Ga, \kx)$ by 
$$
\cb \biggl(\, \sum_{i=m+1}^n u_i \phi(b_i), \, \sum_{j=m+1}^n v_j \phi(b_j) \biggr) = \prod_{m+1\le i < j \le n} \chi(b_i,b_j)^{-u_i v_j}.
$$
Then $\cbsh( \phi(b_i), \phi(b_j) ) =  \chi(b_i,b_j)^{-1}$ for $i,j \in [m+1,n]$, and so
$$
\chi_\cb(b_i,b_j) = \begin{cases}
\chi(b_i,b_j) &\quad (i \le m \; \text{or} \; j \le m)  \\  1 &\quad (i,j > m)
\end{cases}
$$
for $i,j \in [1,n]$.
It follows that
$$\tw_\Ga(A) \subseteq \AD(A^\cb) = \langle \bfq(\chi_\cb) \rangle =
\langle \chi(b_i,b_j) \mid i \le m \; \text{or} \; j \le m \rangle  = \langle \chi(\ker \phi, \, \Zset^n) \rangle,
$$
which completes the proof of the proposition.
\end{proof}

If the homomorphism $\phi$ in \thref{twGa.sandwich} is an isomorphism, equation \eqref{twGa.eqn} says that $\tw_\Ga(A)$ is trivial. This yields the following.

\bco{trivial.tw}
Let $\bfq \in M_n(\kx)$ be a multiplicatively skew-symmetric matrix, and give $\Obfqkn$ and $\TT^n_\bfq(\KK)$ the standard $\Zset^n$-gradings for which the canonical generators $y_i$ are homogeneous of degree $e_i$. For every inclusion $\Obfqkn \subseteq A \subseteq \TT^n_\bfq(\KK)$ of $\Zset^n$-graded $\KK$-algebras, we have
$$
\tw_{\Zset^n}(A) = \langle 1\rangle.
$$
\eco

Many algebras $A$ with a quantum cluster $\Obfqkn \subseteq A \subseteq \TT^n_\bfq(\KK)$ are not graded by $\Zset^n$ but by a proper homomorphic image $\Ga$ of $\Zset^n$. We shall see that $\tw_\Ga(A)$ is typically nontrivial in such cases. In particular, in the case when $\Ga = \Zset$, the invariant $\tw_\Ga(A)$ reduces to $\AD(A)$, as follows.

\bpr{Ga=Z}
Let $\bfq \in M_n(\kx)$ be a multiplicatively skew-symmetric matrix, and suppose $\Obfqkn$ and $\TT^n_\bfq(\KK)$ have $\Zset$-gradings such that the canonical generators $y_i$ are homogeneous. For every inclusion $\Obfqkn \subseteq A \subseteq \TT^n_\bfq(\KK)$ of $\Zset$-graded $\KK$-algebras, we have
$$
\tw_{\Zset}(A) = \langle \bfq\rangle.
$$
\epr

\begin{proof} View the given $\Zset$-gradings on $\Obfqkn$ and $\TT^n_\bfq(\KK)$ as induced from the standard $\Zset^n$-gradings along a suitable homomorphism $\phi : \Zset^n \rightarrow \Zset$. As in the proof of \thref{twGa.sandwich}, 
$$
\tw_\Zset(A) = \bigcap_{\cb \in Z^2(\Zset, \kx)} \langle \bfq(\chi_\cb) \rangle,
$$
where $\chi$ is the alternating bicharacter on $\Zset^n$ such that $\chi(e_i, e_j) = q_{ij}$ for all $i$, $j$ and the alternating bicharacters $\chi_\cb$ are defined as in \eqref{chi_c}.

Observe that any alternating bicharacter $\db$ on $\Zset$ is trivial, since $\db(1,1)=1$ and $\db(s,t) = \db(1,1)^{st}$ for all $s,t \in \Zset$. Hence, $\chi_\cb = \chi$ for all $\cb \in Z^2(\Zset, \kx)$, and thus $\tw_\Zset(A) = \langle \bfq(\chi) \rangle = \langle \bfq \rangle$.
\end{proof}

\bex{qWeyl}
The \emph{multiparameter quantized Weyl algebra} $A := A_n^{Q,P}(\KK)$ is given by generators
 $x_1,y_1,\dots,
x_n,y_n$ satisfying the following relations:
\begin{equation}  \label{AnQPrels}
\begin{aligned}
y_iy_j &=p_{ij}y_jy_i \,, & &(\text{all\ }
i,j),  \\ 
x_ix_j &=q_ip_{ij}x_jx_i \,, &\qquad &(i<j),  \\
x_iy_j &=p_{ji}y_jx_i \,, & &(i<j),  \\
x_iy_j &=q_jp_{ji}y_jx_i \,, & &(i>j),  \\
x_jy_j &=1+q_jy_jx_j +\sum_{l<j} (q_l-1)y_lx_l \,, & &(\text{all\ }j),
\end{aligned}
\end{equation}
where $Q = (q_1,\dots,q_n) \in (\kx)^n$ and $P = (p_{ij}) \in M_n(\kx)$ is a multiplicatively skew-symmetric matrix. There are commuting normal elements $z_k := x_k y_k - y_k x_k \in A$ for $k \in [1,n]$, satisfying
\begin{align*}
z_k y_j &= \begin{cases}
q_j y_j z_k &(j \le k),  \\  y_j z_k &(j>k),
\end{cases}
&z_k x_j &= \begin{cases}
q_j^{-1} x_j z_k &(j \le k),  \\ x_j z_k &(j>k),
\end{cases}
&&\forall\, j,k \in [1,n]
\end{align*}
\cite[\S\S2.8, 3.1]{Jor}. Moreover, $A$ has a quantum cluster of the form
\begin{equation}  \label{qWeyl.clu}
\OO_\bfq(\KK^{2n}) \subseteq A \subseteq \TT_\bfq^{2n}(\KK)
\end{equation}
where $\OO_\bfq(\KK^{2n})$ is generated by $z_1,\dots,z_n,y_1,\dots,y_n$ \cite[\S3.1]{Jor}. 

We may grade $A$ by $\Ga := \Zset^n$, with each $x_k$ homogeneous of degree $e_k$ and each $y_k$ homogeneous of degree $-e_k$. If we assign to $\OO_\bfq(\KK^{2n})$ and $\TT_\bfq^{2n}(\KK)$ the $\Zset^n$-gradings induced from their standard $\Zset^{2n}$-gradings along the homomorphism 
$$
\phi : \Zset^{2n} \longrightarrow \Ga, \;\; \phi(e_j) = \begin{cases}
0 &(j\le n),  \\  -e_{j-n} &(j>n),
\end{cases}
\;\; \forall\, j \in [1,2n],
$$
then \eqref{qWeyl.clu} is an inclusion of $\Ga$-graded $\KK$-algebras. We may write $\OO_\bfq(\KK^{2n}) = \OO_\chi(\KK^{2n})$ and $\TT_\bfq^{2n}(\KK) = \TT_\chi^{2n}(\KK)$ where $\chi$ is the alternating bicharacter on $\Zset^{2n}$ such that 
\begin{align*}
\chi(e_k,e_j) &= 1, &\chi(e_{n+k}, e_{n+j}) &= p_{kj} \,, &\chi(e_k, e_{n+j}) &= \begin{cases}
q_j &(j \le k),  \\  1 &(j>k),
\end{cases}
\end{align*} 
for all $j,k \in [1,n]$.
Since $\ker \phi = \bigoplus_{k=1}^n \Zset e_k$, we conclude from \eqref{twGa.eqn} of \thref{twGa.sandwich} that
\begin{equation}  \label{twGa.qWeyl}
\tw_{\Ga}(A_n^{Q,P}(\KK)) = \langle q_1,\dots,q_n \rangle.
\end{equation}
\eex
\sectionnew{The twist invariant for quantum nilpotent algebras and quantum Schubert cell algebras}
In this section we recall results on the structure of CGL extensions (also known as quantum nilpotent algebras), which form
a large class of algebras that contains as special cases various families of algebras from the theory of quantum groups. 
A theorem is proved for the effective calculation of their twist invariants and is illustrated for the 
algebras of quantum matrices and more generally the quantum Schubert cell algebras for all symmetrizable 
Kac--Moody algebras.

We apply the twist invariant to prove that many quantized Weyl algebras are not uniparameter CGL extensions. This means that the scalars 
in the leading terms of the commutation relations of these CGL extensions do not form a cyclic group.

\subsection{Quantum nilpotent algebras and their AD invariants}
The class of \emph{quantum nilpotent algebras} consists of certain iterated skew polynomial extensions of the form
\begin{equation} 
\label{itOre}
A := \KK[x_1][x_2; \sig_2, \delta_2] \cdots [x_N; \sig_N, \delta_N],
\end{equation}
where the $\sig_k$ are $\KK$-algebra automorphisms and the $\de_k$ are $\KK$-linear $\sig_k$-derivations. The precise conditions were formulated by Launois-Lenagan-Rigal in \cite[Definition 3.1]{LLR} and baptized ``CGL extensions".

\bde{CGL} An iterated skew polynomial extension \eqref{itOre}
is called a \emph{CGL extension} 
 if it is equipped with a rational action of a $\KK$-torus $\HH$ 
by $\KK$-algebra automorphisms satisfying the following conditions:
\begin{enumerate}
\item[(i)] The elements $x_1, \ldots, x_N$ are $\HH$-eigenvectors.
\item[(ii)] For every $k \in [2,N]$, $\de_k$ is a locally nilpotent 
$\sig_k$-derivation of the algebra 
\begin{equation} 
\label{Ak}
A_{k-1} = \KK[x_1][x_2; \sig_2, \delta_2] \cdots [x_{k-1}; \sig_{k-1}, \delta_{k-1}]
\end{equation}
\item[(iii)] For every $k \in [1,N]$, there exists $h_k \in \HH$ such that 
$\sig_k = (h_k \cdot)|_{A_{k-1}}$ and the $h_k$-eigenvalue of $x_k$, to be denoted by $\la_k$, is not a root of unity.
\end{enumerate}

Conditions (i) and (iii) imply that 
$$
\sig_k(x_j) = \la_{kj} x_j \; \; \mbox{for some} \; \la_{kj} \in \kx, \; \; \forall\, 1 \le j < k \le N.
$$
We then set $\la_{kk} :=1$ and $\la_{jk} := \la_{kj}^{-1}$ for $j< k$, obtaining a multiplicatively skew-symmetric 
matrix $\lab := (\la_{kj}) \in M_N(\kx)$.

Denote the character group of the torus $\HH$ by $\xh$ and treat it as an additive abelian group.  The action of 
$\HH$ on $A$ gives rise to an $\xh$-grading of $A$ (e.g., \cite[Lemma II.2.11]{BG}), and the $\HH$-eigenvectors in $A$
are precisely the nonzero homogeneous elements with respect to this grading. 
\ede 

For a function $\eta : [1,N] \to \Zset$, define the predecessor function for its level sets, $p = p_\eta : [1,N] \to [1,N] \sqcup \{ - \infty \}$, by
\[
p(k) := \max \{ j <k \mid \eta(j) = \eta(k) \},
\]
where we set $\max\, \varnothing = - \infty$. For $k \in [1,N]$, denote
\[
O_-(k) := \max \{ l \in \Znn \mid p^l(k) \ne -\infty \}.
\]

An element $y$ of an algebra $A$ is called {\em{prime}} if it is a nonzero nonunit,
\[
Ay = yA, \quad \mbox{and} \quad A/A y \; \; \mbox{is a domain}. 
\] 
We will need the following theorem for the homogeneous prime elements of the algebras 
$A_k$ from \eqref{Ak},
which are CGL extensions for the restriction of the $\HH$-action to them. Set $A_0 := \KK$.

\bth{primeCGL} \cite[Theorems 4.3, 4.6]{GoYa} Let $A$ be a CGL extension as in Definition {\rm\ref{dCGL}}.

{\em{(a)}} Each of the algebras $A_k$ contains a unique {\rm(}up to scalar multiples{\rm)} $\xh$-ho\-mo\-gen\-e\-ous prime element 
that does not belong to $A_{k-1}$. This element will be denoted by $y_k$.

{\em{(b)}} There exists a function $\eta : [1,N] \to \Zset$ such that {\rm(}up to scalar multiples{\rm)} the elements 
$y_1, \ldots, y_N$ are given by  
\[
y_k = 
\begin{cases}
y_{p(k)} x_k - c_k, \;\; \text{for some} \;\; c_k \in A_{k-1}\,,   &\mbox{if} \; \;  p(k) \neq - \infty \\
x_k, & \mbox{if} \; \; p(k) = - \infty.
\end{cases}
\]

{\em{(c)}} The elements $y_1, \ldots, y_N$ define embeddings 
\[
\OO_\bfq(\KK^N) \subseteq A \subset \TT^N_\bfq(\KK),
\]
where the entries $q_{jk}$ of the matrix $\bfq$ are given by
$$
q_{jk} := \prod_{l=0}^{O_-(j)} \prod_{m=0}^{O_-(k)} \la_{p^l(j), p^m(k)} \,.
$$
\eth

It follows from part (b) of the theorem that 
\begin{equation}
\label{de-p}
\delta_k \neq 0 \quad \mbox{if and only if} \quad p(k) = -\infty, 
\end{equation}
see \cite[Theorem 3.6]{GoYa2} for details.

\bpr{AD.cgl} If $A$ is a CGL extension as in Definition {\rm\ref{dCGL}}, then
$\AD(A) = \langle \lab \rangle$.
\epr

\begin{proof} It is easy to see that $\lcor \bfq \rcor = \lcor \lab \rcor$.
Hence, it follows from \thref{primeCGL}(c), Observation \ref{ADobs}(1) and \prref{ADqaff} that $\AD(A) = \langle \lab \rangle$.
\end{proof}

\subsection{The twist invariant for quantum nilpotent algebras} 
Consider a CGL extension $A$ as in \deref{CGL}. For $h \in \HH$ and $\al \in \xh$, we will denote 
by $h^\al \in \kx$ the corresponding eigenvalue. Consider the homomorphism 
\begin{equation}  \label{degreehom}
\pi : \Zset^N \to \xh \quad \mbox{given by} \quad
\pi(e_i) := \deg x_i \in \xh, \;\; \forall\, i \in [1,N].
\end{equation}
Denote the subsets of $[1,N]$ indexing the variables $x_k$ with vanishing and non-vanishing derivations by
\[
V(A) := \{ k \in [1,N] \mid \de_k = 0 \}, \quad N(A) := [1,N] \backslash V(A).
\]

\ble{spanXH} Keep the above notation.

{\rm(a)} The set $\{ \pi (e_i) \mid i \in V(A) \}$ generates $\img \pi$.

{\rm(b)} For each $k \in N(A)$, there exists $b_k \in \ker \pi$ such that $b_k - e_k \in \sum_{j=1}^{k-1} \Zset e_j$. If such elements $b_k$ have been chosen, and if
\begin{equation}  \label{piei.ind}
\text{The elements $\pi(e_j)$ for $j \in V(A)$ are $\Zset$-linearly independent,}
\end{equation}
then $(b_k)_{k \in N(A)}$ is a basis for $\ker \pi$.
\ele

\begin{proof}
(a) For $k \in N(A)$, choose $j_k \in [1,k-1]$ and $(m_{k,1}, \ldots, m_{k,k-1}) \in \Znn^{k-1}$ such that $\de_k(x_j) \ne 0$ 
and the monomial $x_1^{m_{k,1}} \cdots x_{k-1}^{m_{k,k-1}}$ appears in $\de_k(x_{j_k})$ (when expressed in the standard PBW basis of $A$) with a nonzero coefficient. 
Comparing the $\xh$-degrees of the two sides of the relation $x_k x_{j_k} - (h_k \cdot x_{j_k}) x_k = \delta_k(x_{j_k})$, 
we obtain
\[
\pi (e_k) + \pi(e_{j_k}) = \sum_{i=1}^{k-1} m_{k,i} \pi (e_i).
\]
Now part (a) follows by induction on $k \in N(A)$.

(b) Keeping the notation above, we can choose
$$
b_k := e_k + e_{j_k} - \sum_{i=1}^{k-1} m_{k,i} e_i \in \ker \pi \,, \;\; \forall\, k \in N(A).
$$

Now assume just that for all $k \in N(A)$, we have chosen $b_k \in \ker \pi$ with $b_k - e_k \in \sum_{j=1}^{k-1} \Zset e_j$. Regardless of \eqref{piei.ind},  the elements $b_k$ are $\Zset$-linearly independent, and
$$
\Zset^N = \sum_{j \in V(A)} \Zset e_j + \sum_{k \in N(A)} \Zset b_k \,.
$$
The remainder of part (b) follows.
\end{proof}

The condition \eqref{piei.ind} does not always hold. For instance, if $A = \OO_q(\KK^N)$ is the standard quantum affine space with $q$ not a root of unity, then $A$ is a CGL extension with respect to the action of $\HH = \kx$ under which $\al.x_i = \al x_i$ for all $\al \in \HH$ and $i \in [1,N]$. Here $V(A) = [1,N]$ and $\xh \cong \Zset$, so the $\pi(e_i)$ are not $\Zset$-linearly independent if $N>1$.

For any CGL extension $A$, there is a canonical choice of a torus, acting rationally on $A$ in the manner required by \deref{CGL}, for which \eqref{piei.ind} does hold. This torus, which is maximal in a suitable sense, was denoted $\Hmax(A)$ in \cite[Eq. (5.3)]{GoYa}, and may be described as follows. For each $k \in N(A)$, choose $j_k \in [1,k-1]$ and $(m_{k,1}, \ldots, m_{k,k-1}) \in \Znn^{k-1}$ as in the proof of \leref{spanXH}. Then $\Hmax(A)$  is the following subgroup of $(\kx)^N$:
\begin{equation}  \label{Hmax}
\Hmax(A) = \{ (\psi_1,\dots,\psi_N) \in (\kx)^N \mid \psi_k \psi_{j_k} = \prod_{i=1}^{k-1} \psi_i^{m_{k,i}}, \; \forall\, k \in N(A) \}.
\end{equation}
It acts on $A$ so that
$$
(\psi_1,\dots,\psi_N).x_i = \psi_i x_i \,, \;\; \forall\, (\psi_1,\dots,\psi_N) \in \Hmax(A), \; i \in [1,N].
$$
The action of the given torus $\HH$ (as in \deref{CGL}) on $A$ factors through $\Hmax(A)$ via a morphism of algebraic groups, and so $A$ is also a CGL extension with respect to the $\Hmax(A)$-action \cite[\S5.2]{GoYa}. The character group $X(\Hmax(A))$ may be identified with $\Zset^{V(A)}$, and the degree homomorphism $\pi : \Zset^N \rightarrow X(\Hmax(A))$ of \eqref{degreehom} is then determined recursively by
\begin{align*}
\pi(e_j) &= e_j \,,  &&j \in V(A),  \\
\pi(e_k) &= - \pi(e_{j_k}) + \sum_{i=1}^{k-1} m_{k,i} \pi(e_i),  &&k \in N(A).
\end{align*}
In particular, the $\pi(e_j)$ for $j \in V(A)$ are the standard basis vectors for $\Zset^{V(A)}$. 

For an algebra $A$ as in \eqref{itOre}, denote the \emph{interval subalgebra}
$$
A_{[j,k]} := \KK \langle x_i \mid j \le  i \le k \rangle, \; \; \forall\,  j,k \in [1,N].
$$
In particular, $A_{[j,k]} = \KK$ if $j \nleq k$.

\bde{symCGL} We call a CGL extension $A$ as in \deref{CGL} {\em{symmetric}} if the following two conditions hold:
\begin{enumerate}
\item[(i)] For all $1 \leq j < k \leq N$,
$$
\de_k(x_j) \in A_{[j+1, k-1]}.
$$
\item[(ii)] For all $j \in [1,N]$, there exists $h^\sy_j \in \HH$ 
such that 
$$
h^\sy_j \cdot x_k = \la_{kj}^{-1} x_k = \la_{jk} x_k, \; \; \forall\,  
k \in [j+1, N]
$$
and $h^\sy_j \cdot x_j = \la^\sy_j x_j$ for some $\la^\sy_j \in \kx$ which is not a root of unity.
\end{enumerate}
\ede
A symmetric CGL extension $A$ has a presentation as a CGL extension 
with the variables $x_k$ in descending order:
\[
A = \KK[x_N] [x_{N-1}; \sigma^*_{N-1}, \de^*_{N-1}] \cdots [x_1; \sigma^*_1, \de^*_1],
\]
see \cite[Corollary 6.4]{GoYa}. It follows from condition (i) in \deref{CGL} and 
condition (i) in \deref{symCGL} that, for $j <k$, the interval subalgebras $A_{[j,k]}$ also have 
presentations as CGL extensions 
\[
A_{[j,k]} = \KK[x_j][x_{j+1}; \sig_{j+1}, \delta_{j+1}] \cdots [x_k; \sig_k, \delta_k]
\]
for the appropriate restrictions of the maps $\sig_{j+1}, \de_{j+1}, \ldots, \sig_k, \de_k$. 

\bth{twGa.cgl} Let $A$ be a CGL extension as in Definition {\rm\ref{dCGL}}, and assume that \eqref{piei.ind} holds {\rm(}for instance, in case $\HH = \Hmax(A)${\rm)}. Let $\chi'$ be the alternating bicharacter on $\Zset^N$ such that $\chi'(e_i, e_j) = \la_{ij}$ for $i,j \in [1,N]$. 

{\em{(a)}} Choose $b_k \in \Zset^N$ as in Lemma {\rm\ref{lspanXH}(b)}. With the above notation,
$$
\tw_{\xh}(A) = \langle \la_k, \, \chi'(b_k, e_i)  \mid k \in N(A), \; i \in V(A), \; i <k \rangle.
$$

{\em{(b)}} If $A$ is a symmetric CGL extension then 
$$
\tw_{\xh}(A) = \langle \la_k \mid k \in N(A) \rangle.
$$
\eth

\begin{proof}[Proof of Theorem {\rm\ref{ttwGa.cgl}(a)}] 
There is an automorphism $\beta$ of $\Zset^N$ such that
$$
\beta(e_k) = \sum_{l=0}^{O_-(k)} e_{p^l(k)}, \;\; \forall\, k \in [1,N],
$$
and $\chi := \chi' \circ (\beta \times \beta)$ is an alternating bicharacter on $\Zset^n$ such that
$$
\chi(e_k,e_j) = q_{kj}, \;\; \forall\, k,j \in [1,N].
$$
In the notation of \deref{Ochi}, \thref{primeCGL}(c) says that
\begin{equation}  \label{sandwich}
\OO_\chi(\KK^N) \subseteq A \subset \TT^N_\chi(\KK).
\end{equation}
Giving $A$ its $\xh$-grading (from the $\HH$-action) and $\OO_\chi(\KK^N)$ and $\TT_\chi(\KK)$ the $\xh$-gradings induced from their canonical $\Zset^N$-gradings along $\phi := \pi \beta$, we see that \eqref{sandwich} is an inclusion of $\xh$-graded $\KK$-algebras. Thus, \thref{twGa.sandwich} implies that
$$
\tw_{\xh}(A) = \langle \chi(\ker \phi, \, \Zset^N) \rangle.
$$
Since $\beta$ is an isomorphism, this yields
\begin{equation}  
\label{twGaCGL}
\tw_{\xh}(A) = \langle  \chi' (\ker \pi, \, \Zset^N) \rangle = \langle \chi'(b_k, e_j) \mid k \in N(A), \; j \in [1,N] \rangle,
\end{equation}
taking account of \leref{spanXH}(b).

Conditions (i) and (iii) in \deref{CGL} imply the following important property of $\chi'$:
\begin{equation}
\label{h-to-chi}
\chi'(e_j, a) = h_j^{\pi(a)}, \; \; \forall\, j \in [2,N], \; a \in \bigoplus_{i < j} \Zset e_i \,.
\end{equation}
It follows from this property that
\begin{equation}
\label{la-k}
\chi'(b_k, e_k) = \chi'(e_k, e_k -b_k) = h_k^{\pi(e_k -b_k)} =
h_k^{\pi(e_k)} = \la_k \,, \;\; \forall\; k \in N(A),
\end{equation}
and that
\begin{equation}
\label{eq-to1}
\chi'(b_k, e_j) = h_j^{- \pi(b_k)} = h_j^0 = 1, \; \; \forall\, k \in N(A), \; j \in [k+1,N].
\end{equation}

For $k \in N(A)$, write
\[
b_k = e_k + a_k \quad \mbox{with} \quad a_k \in \bigoplus_{i=1}^{k-1} \Zset e_i.
\]  
For $k, j \in N(A)$ with $k >j$, we have,
\begin{align*}
\chi'(b_k, e_j) &= \chi'(b_k, b_j - a_j) = \chi'(b_k, b_j) \chi'(b_k, a_j)^{-1} = \chi'(e_k, b_j) \chi'(a_k, b_j) \chi'(b_k, a_j)^{-1} \\
& = \chi'(b_j, a_k)^{-1} \chi'(b_k, a_j)^{-1} \in \langle \la_j, \chi'(b_k, e_s), \chi'(b_j, e_t), \, s, t < j\rangle,
\end{align*}
where in the last equality we used \eqref{la-k}--\eqref{eq-to1}. An inductive argument yields that
for $k, j \in N(A)$ with $k >j$
\[
\chi'(b_k, e_j) \in \langle \la_l, \, \chi'(b_l, e_i)  \mid l \in N(A), \; i \in V(A), \; i < l \leq k \rangle.
\]
Part (a) of the theorem follows from this property and Eqs. \eqref{twGaCGL}, \eqref{la-k}--\eqref{eq-to1}.
\end{proof}
\subsection{Proof of \thref{twGa.cgl}(b)} 
Before we proceed with the proof of part (b) of \thref{twGa.cgl} we establish a lemma. Continue with the previous notation.

\ble{add-ker-elem} Let $A$ be a symmetric CGL extension. For each $k \in N(A)$, there exists $a_k \in \bigoplus_{i=p(k)+1}^{k-1} \Zset e_i$ 
such that 
\begin{align}
\label{bk}
b_k &:= e_{p(k)} - a_k + e_k \in \ker \pi  \\
\label{chi'bkej}
\chi'(b_k,e_j) &= \begin{cases}
1, &\quad \forall\, j \in [1,N] \setminus \{p(k),k\},  \\
\la_k, &\quad j = k,  \\
\la_k^{-1}, &\quad j = p(k).
\end{cases}
\end{align}
\ele

\begin{proof} Fix $k \in N(A)$. Then $p(k) \ne -\infty$ and the interval subalgebra $A' := A_{[p(k),k]}$ is a CGL extension (with the restricted action of $\HH$) of the form
$$
A' = \KK[x'_1] [x'_2; \sig'_2, \de'_2] \cdots [x'_n; \sig'_n, \de'_n]
$$
where $n := k-p(k)+1$ and $x'_i := x_{p(k)+i-1}$ for $i \in [1,n]$, while $\sig'_i$ and $\de'_i$ are the restrictions of $\sig_{p(k)+i-1}$ and $\de_{p(k)+i-1}$ to $A_{[p(k), p(k)+i-2]}$ for $i \in [2,n]$. The $\eta$-function for $A'$ given by \thref{primeCGL}(b), call it $\eta'$, can be chosen as follows:
$$
\eta'(i) = \eta(p(k)+i-1), \;\; \forall\, i \in [1,n],
$$
 see \cite[Corollary 5.6(b)]{GoYa2}. If $p'$ is the corresponding predecessor function, then $p'(n) = 1$. Hence, up to a scalar multiple the final $\xh$-homogeneous prime element of $A'$ in \thref{primeCGL}(a)(b) is
\begin{equation}  \label{ynprime}
y'_n = x'_1 x'_n - c'_n
\end{equation}
where $c'_n \in A_{[p(k), k-1]}$ and all terms in \eqref{ynprime} are homogeneous of the same $\xh$-degree. In fact, $c'_n$ is a scalar multiple of $\de'_n(x'_1) = \de_k(x_{p(k)})$ \cite[Proposition 4.7(b)]{GoYa}, which means that
$$
c'_n \in A_{[p(k)+1, k-1]} \,,
$$
due to the symmetry assumption on $A$.
Moreover, by \cite[Corollary 4.8]{GoYa},
\begin{equation}  \label{yn'comm}
y'_n x_j = (\la_{jk} \la_{j,p(k)})^{-1} x_j y'_n \,, \;\; \forall\, j \in [p(k),k].
\end{equation}

The monomials 
\[
x_{p(k)}^{n_{p(k)}} \ldots x_k^{n_k}, \; \; (n_{p(k)}, \ldots, n_k) \in \Znn^n
\]
form a $\KK$-basis of $A'$, and  the vector of exponents for the leading term of $y'_n$ with respect to the right-to-left lexicographic 
order on $\Znn^n$ is $(1, 0, \ldots, 0, 1)$. 
Let $(0, m_{p(k)+1}, \ldots, m_{k-1}, 0)$ be the vector of exponents for the leading term of $c'_n$, and set
\[
a_k := m_{p(k)+1} e_{p(k)+1} + \cdots + m_{k-1} e_{k-1}.
\]
It follows from the $\xh$-homogeneity of \eqref{ynprime} that \eqref{bk} holds, and from \eqref{la-k}, \eqref{eq-to1} we obtain \eqref{chi'bkej} for $j \in [k,N]$.

As far as the CGL extension $A_{[p(k)+1,k]}$ is concerned, the index for $x_k$ has no predecessor, and so $x_k$ is a homogeneous normal element of that algebra. Similarly, $x_{p(k)}$ is a homogeneous normal element of the algebra $A_{[p(k),k-1]}$. Hence,
$$
x_{p(k)} x_k x_j = \la_{p(k),j} \la_{kj} x_j x_{p(k)} x_k \,, \;\; \forall\, j \in [p(k)+1, k-1].
$$
Taking \eqref{yn'comm} into account, we see that
$$
c'_n x_j = \la_{p(k),j} \la_{kj} x_j c'_n \,, \;\; \forall\, j \in [p(k)+1, k-1].
$$
This equation, together with the relation \cite[Eq. (4.12)]{GoYa} for leading terms of products of monomials in CGL extensions, implies
$$
\prod_{i=j+1}^{k-1} \la_{ij}^{m_i} = \la_{p(k),j} \la_{kj} \prod_{i=p(k)+1}^{j-1} \la_{ji}^{m_i} \,, \;\; \forall\, j \in [p(k)+1, k-1],
$$
which yields
$$
\chi'( b_k, e_j ) = 1, \,, \;\; \forall\, j \in [p(k)+1, k-1].
$$

Condition (ii) of \deref{symCGL} implies the following analog of \eqref{h-to-chi}:
$$
\chi'(e_j, a) = (h^*_j)^{\pi(a)} \,, \;\; \forall\, j \in [1, p(k)-1], \;\; a \in \bigoplus_{i>j} \Zset e_i \,.
$$
As with \eqref{la-k} and \eqref{eq-to1}, it follows that
\begin{align*}
\chi'(b_k, e_{p(k)}) &= \la^*_{p(k)} = \la_k^{-1} \,,  \\
\chi'(b_k, e_j) &= 1, \qquad\qquad\qquad \forall\, j \in [1, p(k)-1],
\end{align*}
where in the first equality we used \cite[Proposition 5.8]{GoYa2}.
This completes the proof of the lemma. 
\end{proof}

\begin{proof}[Proof of Theorem {\rm\ref{ttwGa.cgl}(b)}] For $k \in N(A)$, let $b_k := e_{p(k)} - a_k + e_k$ as in \eqref{bk}, and note that these elements satisfy the hypotheses of \leref{spanXH}(b). In view of \eqref{chi'bkej},
$$
\chi'(b_k, e_j) \in \langle \la_k \rangle, \;\; \forall\, k \in N(A), \;\; j \in [1,N].
$$
Part (b) of the theorem follows from this and part (a). 
\end{proof}

\bex{OlpMn2}
It is well known that $A := \OlpMn$ can be expressed as a CGL extension provided $\la$ is a non-root of unity (cf.~\cite[\S5.6]{GLet} and \cite[Corollary 3.8]{LLR}) with the generators $X_{ij}$ adjoined in lexicographic order and the torus $\HH := (\kx)^{2n}$ acting so that
$$
\al.X_{ij} = \al_i \al_{n+j} X_{ij} \,, \;\; \forall\, \al = (\al_1,\dots,\al_{2n}) \in \HH, \;\; i,j \in [1,n].
$$
Equally, $A$ can be expressed as a CGL extension with the generators adjoined in reverse lexicographic order and the same $\HH$-action, which shows that this algebra is a symmetric CGL extension. We may identify $\xh$ with $\Zset^{2n}$, and we note that the $\xh$-grading of $A$ equals the $\Zset^{2n}$-grading given in \exref{OlpMn}. Observe that
$$
V(A) = [1,n] \cup \{ in+1 \mid i \in [1,n-1] \},
$$
and that the degree homomorphism $\pi : \Zset^{n^2} \rightarrow \xh$ satisfies $\pi(j) = e_1+ e_{n+j}$ for $j \in [1,n]$ while $\pi(in+1) = e_{i+1} + e_{n+1}$ for $i \in [1,n-1]$. Hence, \eqref{piei.ind} holds.

The calculations in \cite[\S5.6]{GLet} show that $\la_k = \la$ for all $k \in N(A)$. Therefore, by \prref{AD.cgl} and \thref{twGa.cgl},
\begin{equation}  \label{ADtw.OlpMn}
\AD(\OlpMn) = \langle \la, \bfp \rangle \qquad\text{and}\qquad \tw_{\xh}(\OlpMn) = \langle \la \rangle,
\end{equation}
assuming $\la$ is a non-root of unity.
Specializing to the standard quantum matrix algebra yields
\begin{equation}  \label{ADtw.OqMn}
\AD(\OO_q(M_n(\KK))) = \langle q \rangle \qquad\text{and}\qquad \tw_{\xh}(\OO_q(M_n(\KK))) = \langle q^2 \rangle
\end{equation}
provided $q$ is not a root of unity.
\eex

\subsection{Uniparameter CGL extensions}
\bde{uni} We will say that a CGL extension $A$ as in \deref{CGL} is \emph{uniparameter} if there is some $q\in \kx$ such that all $\la_{ij} \in \langle q \rangle$, 
that is, $\langle \lab \rangle$ is a cyclic group. 
\ede

For example, the standard quantum matrix algebras $\OO_q(M_n(\KK))$ are uniparameter, and all the multiparameter quantum matrix algebras $\OlpMn$ are cocycle twists of uniparameter ones, as shown in \exref{OlpMn}. On the other hand, many quantum algebras cannot be so reduced to uniparameter algebras, as we now illustrate.

\bex{qWeyl2}
Let $A := A_n^{Q,P}(\KK)$ as in \exref{qWeyl}, and assume that all the $q_k$ are non-roots of unity. As is well known, $A$ is then a CGL extension with variables in the order $x_1,y_1,\dots, x_n,y_n$, equipped with the action of $\HH := (\kx)^n$ by algebra automorphisms such that
\begin{equation}  \label{Hact.qWeyl}
(\al_1,\dots,\al_n) . x_i = \al_i x_i \; \; \text{and} \; \; (\al_1,\dots,\al_n) . y_i = \al_i^{-1} y_i \,,  \; \; \forall\, i \in [1,n]
\end{equation}
(cf.~\cite[\S5.10]{GLet}, \cite[Corollary 3.8]{LLR}). We can identify the $\xh$-grading of $A$ with the $\Zset^n$-grading of \exref{qWeyl}, and so by \eqref{twGa.qWeyl},
$$
\tw_{\xh}(A_n^{Q,P}(\KK)) = \langle q_1,\dots,q_n \rangle.
$$
(This can also be obtained from \thref{twGa.cgl} after rewriting $A$ as a symmetric CGL extension, which can be done by adjoining the generators in the order $y_n,\dots,y_1, x_1,\dots,x_n$.) It follows that whenever the group $\langle q_1,\dots,q_n \rangle$ fails to be cyclic, $A_n^{Q,P}(\KK)$ is not cocycle twist equivalent to a uniparameter CGL extension.
\eex

\subsection{Quantum Schubert cell algebras}
Let $\g$ be a symmetrizable Kac--Moody algebra with set of simple roots $\al_1, \ldots, \al_r$, 
Weyl group $W$, and root lattice $Q$. Denote by $U_q(\g)$ the 
corresponding quantized universal enveloping algebra over $\KK$, for a 
deformation parameter $q \in \kx$ which is not a root of unity. Let 
$\{K_i^{\pm 1}, E_i, F_i \mid i \in [1,r] \}$ be the generators of $U_q(\g)$. Denote by 
\[
\lcor ., . \rcor : Q \times Q \longrightarrow \Zset
\]
the invariant symmetric form, normalized by $\| \al \|^2 := \lcor \al, \al \rcor = 2$ for a short root $\al$. 

To each element $w \in W$, one associates a quantum Schubert cell subalgebra $U^q_-[w] \subset U_q(\g)$, \cite[\S 39.3]{L}.
Given a reduced expression $w = s_{i_1} \ldots s_{i_l}$, define the roots
\[
\beta_k := s_{i_1} \ldots s_{i_{k-1}} (\al_{i_k}), \quad \forall\,  k \in [1, l]
\]
of $\g$ and the root vectors 
\[
F_{\be_k} := T_{i_1} \ldots T_{i_{k-1}} (E_{i_k}), 
\quad \forall\, k \in [1,l]
\]
of $U_q(\g)$. Here $T$ denotes Lusztig's action of the braid group of $\g$ on $U_q(\g)$, \cite[Ch. 5]{L}. The algebra $U^q_-[w]$ is defined as
the subalgebra of $U_q(\g)$ generated by $F_{\be_1}, \ldots, F_{\be_l}$. This subalgebra is independent of the choice of 
a reduced expression of $w$. The algebras of quantum matrices $\OO_q(M_n(\KK))$ are recovered as special cases 
for $\g = \sl_{2n}$ and a particular choice of $w \in S_{2n}$, \cite[Sect. 4]{Y}.

The algebra $U_q(\g)$ is $Q$-graded by 
\[
\deg F_{\be_k} := \be_k \,, \quad \forall\, k \in [1,l].
\]
The grading gives rise to a rational action of a torus $\HH := (\kx)^r$ by algebra automorphisms obtained by identifying
$\xh \cong Q$. The algebra $U^q_-[w]$ is a symmetric CGL extension when its generators are adjoined 
in the order $F_{\be_1}, \ldots, F_{\be_l}$, \cite[Lemma 9.1]{GoYa2}. 
The corresponding scalars $\la_1, \ldots, \la_l$ are given by
\begin{equation}
\label{scalar-la}
\la_k = q^{ - \| \be_k \|^2} = q^{ - \| \al_{i_k} \|^2}, \quad \forall\, k \in [1,l],
\end{equation}
see \cite[Eq. (9.15)]{GoYa2}. These facts were stated in \cite{GoYa2} in the case where $\g$ is finite dimensional, but the 
proofs in the symmetrizable Kac--Moody case are identical.

Denote the support of $w \in W$ by
\[
\supp (w) := \{ j \in [1,r ] \mid j = i_k \; \; \mbox{for some} \; \; k \in [1,l] \},
\]
and set
\[
d(w) : = {\mathrm{gcd}}( \| \al_j \|^2 , \; j \in \supp (w) ). 
\]

\thref{twGa.cgl}(b) and \eqref{scalar-la} imply at once the following result.

\bpr{twGamUqw} For all symmetrizable Kac--Moody algebras $\g$ with root lattice $Q$, Weyl group elements $w$, and non-roots 
of unity $q \in \kx$, 
\[
\tw_Q(U^q_-[w]) = \langle q^{ d(w)} \rangle.
\]
\epr
\sectionnew{General twist invariants}

For a $\Zset$-graded algebra $A$, the invariant $\tw_\Zset(A)$ may not distinguish $A$ from twists of $A$ by $2$-cocycles on finer grading groups, as illustrated by \prref{Ga=Z}. In order to recover cocycle twist invariant information within a $\Zset$-graded invariant, we work with gradings of $A$ by character groups of maximal tori of suitable automorphism groups of $A$ with algebraic group structures. For this process to function, the homogeneous components of the $\Zset$-grading on $A$ need to be finite dimensional. We can enlarge the setting to include filtered algebras, and we begin in that setting.

\subsection{Filtered automorphism groups}  \label{filtautgroups}
\bde{defAutfl} All filtered $\KK$-algebras $A$ that we consider will be nonnegatively filtered, and the filtration will be denoted $F_\bullet A = (F_n A)_{n\ge0}$. The \emph{filtered automorphism group} of such a filtered algebra is
$$
\autfl A := \{ \phi \in \Aut_\KK A \mid \phi(F_n A) = F_n A, \;\; \forall\, n \ge 0 \},
$$
where $\Aut_\KK A$ denotes the group of all $\KK$-algebra automorphisms of $A$. Recall that the filtration $F_\bullet A$ is called \emph{locally finite} if each of the subspaces $F_n A$ is finite dimensional over $\KK$. 
\ede

\ble{rVW} 
Let $A$ be an affine $\KK$-algebra and $V \subseteq W$ finite dimensional $\KK$-subspaces of $A$ such that $V$ generates $A$ as an algebra. Set
$$
\AA_V := \{ \phi \in \Aut_\KK A \mid \phi(V) = V \} \qquad\text{and}\qquad \AA_W := \{ \phi \in \Aut_\KK A \mid \phi(W) = W \},
$$
and let
$$
\rho_V : \AA_V \longrightarrow \GL(V) \qquad\text{and}\qquad \rho_W : \AA_W \longrightarrow \GL(W)
$$
be the homomorphisms given by restriction to $V$ and $W$, respectively.

{\rm(a)} $\rho_V(\AA_V)$ and $\rho_V(\AA_V \cap \AA_W)$ are {\rm(}Zariski-{\rm)}closed subgroups of $\GL(V)$, while $\rho_W(\AA_V)$ and $\rho_W(\AA_V \cap \AA_W)$ are closed subgroups of $\GL(W)$.

{\rm(b)} The map $\rho_{WV} : \rho_W(\AA_V \cap \AA_W) \rightarrow \rho_V(\AA_V \cap \AA_W)$ given by restriction to $V$ is an isomorphism of algebraic groups.
\ele

\begin{proof}  (a) Let $\Ahat := T(V)$ be the tensor algebra of $V$, graded by $\Znn$ in the standard way, so that $\Ahat_m = V^{\otimes m}$, and let $\pi : \Ahat \rightarrow A$ be the unique $\KK$-algebra homomorphism which is the identity on $V$. For $g\in \GL(V)$, let $\ghat$ be the unique $\KK$-algebra automorphism of $\Ahat$ which restricts to $g$ on $V$. Observe that $g \in \rho_V(\AA_V)$ if and only if $\ghat$ induces a $\KK$-algebra automorphism of $A$, if and only if $\ghat(\ker\pi) = \ker\pi$. The last condition is equivalent to $\ghat(K_m) = K_m$ for all $m\ge0$, where $K_m := \Ahat_{\le m} \cap \ker\pi$.

Now $H_m := \{ h \in \GL(\Ahat_{\le m}) \mid h(K_m) = K_m\}$ is a closed subgroup of $\GL(\Ahat_{\le m})$, and the set
$$G_m := \{ g \in \GL(V) \mid \ghat(K_m) = K_m \}$$
is the inverse image of $H_m$ under the morphism 
$$
\tau_m : \GL(V) \longrightarrow \GL(\Ahat_{\le m}), \quad 
g \, \longmapsto \, \bigoplus_{n=0}^m \, (g^{\otimes n})|_{V^{\otimes n}} = \ghat|_{\Ahat_{\le m}},
$$
so $G_m$ is closed in $\GL(V)$. Therefore $\rho_V(\AA_V) = \bigcap_{m\ge0} G_m$ is closed in $\GL(V)$.

Similarly, $\rho_W(\AA_W)$ is closed in $\GL(W)$. Then also
$$
\rho_W(\AA_V \cap \AA_W) = \{ f \in \rho_W(\AA_W) \mid f(V) = V \}
$$
is closed in $\GL(W)$. Choose $m_0 \in \Znn$ large enough so that $W \subseteq V^{m_0}$, and observe that
\begin{multline*}
\rho_V(\AA_V \cap \AA_W) =  \\
\bigcap_{m \ge m_0} \tau_m^{-1} \bigl( \{ l \in \GL(\Ahat_{\le m}) \mid l(K_m) = K_m \;\; \text{and} \;\; l(\pi^{-1}(W)) = \pi^{-1}(W) \} \bigr).
\end{multline*}
It follows that $\rho_V(\AA_V \cap \AA_W)$ is closed in $\GL(V)$.

(b) Injectivity of $\rho_{WV}$ is clear, since automorphisms of $A$ are determined by their actions on $V$. Set
$$
E_V := \rho_V(\AA_V \cap \AA_W) \qquad\text{and}\qquad E_W := \rho_W(\AA_V \cap \AA_W),
$$
and choose $m \in \Znn$ such that $W \subseteq V^m$. Let $\mu: V^{\otimes m} \rightarrow V^m$ be the multiplication map. Set
$$
C := \{ f \in \GL(V^{\otimes m}) \mid f(\ker\mu) = \ker\mu \} \quad\text{and}\quad D := \{ f\in \GL(V^m) \mid f(W) = W \},
$$
which are closed subgroups of $\GL(V^{\otimes m})$ and $\GL(V^m)$, respectively. Restriction from $V^m$ to $W$ provides an algebraic group morphism $r : D \rightarrow \GL(W)$, and induction from $V^{\otimes m}$ to $V^m$ along $\mu$ provides an algebraic group morphism $s : C \rightarrow \GL(V^m)$. Finally, we have a morphism $t : \GL(V) \rightarrow \GL(V^{\otimes m})$ given by $t(g) = g^{\otimes m}$.

If $\phi\in \AA_V \cap \AA_W$ and $g = \rho_V(\phi) = \phi|_V$, then since
$$
\mu t(g)(v_1 \otimes \cdots \otimes v_m) = g(v_1) g(v_2) \cdots g(v_m) = \phi(v_1 v_2 \cdots v_m)
$$
for all $v_1,\dots,v_m \in V$, we see that $\mu t(g) = \phi \mu$. In particular, it follows that $t(g) \in C$. Hence, $s$ composed with $t|_{E_V}$ provides a morphism $st : E_V \rightarrow \GL(V^m)$. For $g$ and $\phi$ as above, $st(g) \mu = \mu t(g) = \phi \mu$, whence $st(g) = \phi|_{V^m}$. Since $\phi \in \AA_W$, it follows that $st(g) \in D$ and $rst(g) = \phi|_W = \rho_W(\phi)$, which means that $rst(g) \in E_W$ and $\rho_{WV} rst(g) = \phi|_V = g$. Therefore $\rho_{WV}$ is bijective and $\rho_{WV}^{-1} = rst|_{E_V}$, which is a morphism. This establishes part (b).
\end{proof}

\bpr{autfl.alg.grp} 
Let $A$ be an affine locally finite filtered $\KK$-algebra, and for $d\in\Znn$ let $\rho_d : \autfl A \rightarrow \GL(F_d A)$ be the homomorphism given by restriction to $F_d A$.

{\rm(a)} If $F_d A$ generates $A$ {\rm(}as an algebra\/{\rm)}, then $\rho_d$ is a group isomorphism from $\autfl A$ onto a {\rm(}Zariski-{\rm)}closed subgroup of $\GL(F_d A)$.

{\rm(b)} If $F_d A$ generates $A$ and $e\ge d$, then the map $\rho_{ed} : \img\rho_e \rightarrow \img\rho_d$ given by restriction to $F_d A$ is an isomorphism of algebraic groups.

{\rm(c)} Choose $d\in \Znn$ such that $F_d A$ generates $A$, and make $\autfl A$ into an affine algebraic group via the isomorphism $\rho_d$. Then the action of $\autfl A$ on $A$ is a rational action.
\epr

\begin{proof} Fix $d \le e$ in $\Znn$ such that $F_dA$ generates $A$. Set $V := F_d A$ and $W := F_eA$. For $n \in \Znn$, define
$$
\AA_n := \{ \phi \in \Aut_\KK A \mid \phi(F_n A) = F_n A \},
$$
and let $\rho'_n : \AA_n \rightarrow \GL(F_nA)$ be the homomorphism given by restriction to $F_nA$.

(a) By \leref{rVW}(a), $\rho'_d(\AA_d)$ is a closed subgroup of $\GL(V)$. It is clear that $\rho'_d$ and $\rho_d$ are injective. Now $\img \rho_d = \rho'_d(\autfl A)$, so it follows from the injectivity of $\rho'_d$ that
$$
\img \rho_d = \bigcap_{n=0}^\infty \rho'_d(\AA_d \cap \AA_n).
$$
For $n \ge d$, we have $\rho'_d(\AA_d \cap \AA_n)$ closed in $\GL(V)$ by \leref{rVW}(a), while for $n<d$, closedness holds because
$$
\rho'_d(\AA_d \cap \AA_n) = \{ f \in \rho'_d(\AA_d) \mid f(F_nA) = F_nA \}.
$$
Therefore $\img \rho_d$ is closed in $\GL(V)$.

(b) By \leref{rVW}(b), the map $\rho'_{ed} : \rho'_e(\AA_d \cap \AA_e) \rightarrow \rho'_d(\AA_d \cap \AA_e)$ given by restriction to $V$ is an isomorphism of algebraic groups. Since $\rho_{ed}$ is the restriction of $\rho'_{ed}$ to $\img \rho_e$, it follows that $\rho_{ed}$ is an algebraic group isomorphism of $\img \rho_e$ onto a closed subgroup of $\rho'_d(\AA_d \cap \AA_e)$. But clearly $\rho_{ed}(\img \rho_e) = \img \rho_d$.

{\rm(c)} For each $e\ge d$, the space $F_e A$ is an $\autfl A$-stable finite dimensional subspace of $A$, and the map $\rho_e : \autfl A \rightarrow \GL(F_e A)$ is a morphism by part (b).
\end{proof}

We now impose the assumption that the field $\KK$ is infinite, so that rational actions of $\KK$-tori on $\KK$-algebras correspond to $\xh$-gradings (e.g., \cite[Lemma II.2.11]{BG}).

\bde{deftwgen} Suppose $A$ is an affine locally finite filtered $\KK$-algebra. Choose $d \in \Znn$ such that $F_d A$ generates $A$, and give $\autfl A$ the structure of an affine algebraic group as in \prref{autfl.alg.grp}(a). By parts (b),(c) of the proposition, this structure is independent of the choice of $d$, and the action of $\autfl A$ on $A$ is rational.
\subsection{General twist invariants of filtered algebras}
Let $\HH$ be a maximal torus of $\autfl A$. The restricted action of $\HH$ on $A$ is rational, so $A$ receives a corresponding $\xh$-grading.

Now suppose $A$ is also a  semiprime right Goldie ring. The \emph{general twist invariant} of $A$ is
$$\tw(A) = \tw(A,F_\bullet) := \tw_{\xh}(A),$$ 
which is well defined by following lemma
\ede

\ble{tw.indepH} Let $A$ be an affine, locally finite, filtered, semiprime right Goldie $\KK$-algebra. Then $\tw(A)$ is independent of the choice of a maximal torus $\HH$ of $\autfl A$.
\ele

\begin{proof} Any two maximal tori $\HH_1$ and $\HH_2$ of $\autfl A$ are conjugate, as follows by applying \cite[Corollary 11.3]{Bor} to the connected group $(\autfl A)^\circ$, so $\HH_2 = \phi \HH_1 \phi^{-1}$ for some $\phi \in \autfl A$. The isomorphism $\theta : \HH_1 \rightarrow \HH_2$ given by conjugation by $\phi$ induces an isomorphism $\theta^* : \Ga_2 := X(\HH_2) \rightarrow \Ga_1 := X(\HH_1)$. We check that $\phi$ provides an isomorphism from ($A$ with $\Ga_1$-grading) to ($A$ with $\Ga_2$-grading), in the sense that for any $d\in \Ga_2$, the automorphism $\phi$ restricts to a vector space isomorphism of $A_{\theta^*(d)}$ onto $A_d$. Namely, if $a \in A_{\theta^*(d)}$, then
$$
\theta(h_1).\phi(a) = \phi h_1 \phi^{-1} . \phi(a) = \phi(h_1.a) = \phi \bigl( \bigl[ \theta^*(d)(h_1) \bigr] a \bigr) = \bigl[ d(\theta(h_1)) \bigr] \phi(a), \;\; \forall\, h_1 \in \HH_1 \,,
$$
whence $\phi(a) \in A_d$. By symmetry, $\phi^{-1}$ maps $A_d$ to $A_{\theta^*(d)}$, establishing the claim. It follows that $\tw_{\Ga_1}(A) = \tw_{\Ga_2}(A)$.
\end{proof}  

\bex{tw.qaff} Consider $A := \Obfqkn$, where $\bfq \in M_n(\kx)$ is an arbitrary multiplicatively skew-symmetric matrix, and give $A$ the standard filtration in which $F_kA$ is the $\KK$-span of the monomials $y_1^{m_1} \cdots y_n^{m_n}$ with $m_1 + \cdots + m_n \le k$. Then $\autfl A$ becomes an algebraic group isomorphic, by restriction, to a closed subgroup of $\GL(F_1A)$. Identify the torus $\HH := (\kx)^n$ in the standard way with a subgroup of $\autfl A$, where
$$
(\al_1, \dots, \al_n)(y_i) = \al_i y_i \,, \;\; \forall\,  (\al_1, \dots, \al_n) \in \HH, \;\; i \in [1,n].
$$
it is easily checked that $\HH$ is a maximal torus of $\autfl A$. When $\xh$ is identified with $\Zset^n$ in the canonical manner, the corresponding grading of $A$ is the one for which $\deg y_i = e_i$ for all $i$. Thus, we conclude from \exref{twGam.qaff} that
$$
\tw(A) = \tw_{\Zset^n}(A) = \{1\}.
$$
\eex

\bth{invar-tw} Let $A$ be an affine, locally finite, filtered, semiprime right Goldie $\KK$-algebra 
and $\HH \subseteq \autfl A$  a $\KK$-torus acting rationally on $A$ by algebra automorphisms. For every $2$-cocycle
$\cb \in Z^2(X(\HH),\kx)$,
\[
\tw(A) = \tw(A^\cb).
\]
\eth

Here $\tw(A^\cb)$ is defined because $F_\bullet A$ is also a filtration of $A^\cb$, and $\HH$ is a subtorus of $\autfl A^\cb$.

\begin{proof} The subtorus $\HH \subseteq \autfl A$ is contained in a maximal torus $\wt{\HH} \subseteq \autfl A$. Using the 
projection $\pi : X(\wt{\HH}) \to \xh$ define the $2$-cocycle
\[
\wt{\cb} := \cb \circ (\pi \times \pi) \in Z^2( X(\wt{\HH}),\kx).
\]
Clearly $A^{\wt{\cb}} \cong A^\cb$ and thus, 
\[
\tw(A^\cb) = \tw(A^{\wt{\cb}}) = \tw_{X(\wt{\HH})}(A^{\wt{\cb}}) = \tw_{X(\wt{\HH})}(A)= \tw(A).  \qedhere
\]
\end{proof}
\subsection{General twist invariants of graded algebras}
\begin{observations} \label{autfl-autgr}
Recall that nonnegatively $\Zset$-graded algebras are turned into nonnegatively filtered algebras in a canonical way. Namely, if $A = \bigoplus_{n=0}^\infty A_n$ is a $\Znn$-graded $\KK$-algebra, we equip $A$ with the associated filtration $F_\bullet A$ where $F_nA := \bigoplus_{k=0}^n A_k$ for all $n \ge 0$. Then $\autfl A$ contains as a subgroup the group
$$
\autgr A := \{ \phi \in \Aut_\KK A \mid \phi(A_n) = A_n \,, \;\; \forall\, n \ge 0 \}
$$
of graded $\KK$-algebra automorphisms of $A$.

Now suppose $A$ is affine and its grading is \emph{locally finite}, meaning that each $A_n$ is finite dimensional over $\KK$. If we give $\autfl A$ the structure of an affine algebraic group via \prref{autfl.alg.grp}, then $\autgr A$ is a Zariski-closed subgroup. More precisely, if $d \in \Znn$ is chosen so that $F_dA$ generates $A$ as an algebra, then restriction to $F_dA$ provides an isomorphism $\rho_d$ of $\autfl A$ onto a closed subgroup of $\GL(F_dA)$, and $\rho_d$ restricts to an isomorphism of $\autgr A$ onto the following closed subgroup of $\img \rho_d$:
$$
\{ f \in \img \rho_d \mid f(A_k) = A_k \,, \;\; \forall\, k \le d \} \cap \bigcap_{e>d} \rho_{ed}^{-1} \bigl( \{ g \in \img \rho_e \mid g(A_e) = A_e \} \bigr),
$$
in the notation of \prref{autfl.alg.grp}. Therefore $\autgr A$ becomes an affine algebraic group.

Thus, if $A$ is an affine locally finite $\Znn$-graded $\KK$-algebra and also a semiprime right Goldie ring, $\tw(A)$ is defined. In setting up $\tw(A)$, one  may use a maximal torus of $\autgr A$ since any such is also a maximal torus of $\autfl A$, as we now show.
\end{observations}

\ble{maxtor.autgr}
Let $A$ be an affine locally finite $\Znn$-graded $\KK$-algebra. Give $\autfl A$ and $\autgr A$  affine algebraic group structures as in Proposition {\rm\ref{pautfl.alg.grp}} and Observations {\rm\ref{autfl-autgr}}. Then any maximal torus of $\autgr A$ is also a maximal torus of $\autfl A$.
\ele

\begin{proof} Let $\HH$ be a maximal torus of $\autgr A$ and $\HH'$ a subtorus of $\autfl A$ which contains $\HH$. We must show $\HH' = \HH$.

Since the filtration on $A$ is the one induced from the grading, the associated graded algebra of $(A,F_\bullet)$ is canonically isomorphic to $A$ as a graded algebra. Consequently, automorphisms of $(A, F_\bullet)$ canonically induce graded automorphisms of $A$, yielding a canonical homomorphism
$$
\ga : \autfl A \longrightarrow \autgr A.
$$
Observe that $\ga$ is a morphism of algebraic groups, and that $\ga$ restricts to the identity on $\autgr A$.

Now $\ga(\HH')$ is a subtorus of $\autgr A$ containing $\ga(\HH) = \HH$, so $\ga(\HH') = \HH$. On the other hand, if we identify $\autfl A$ with the closed subgroup $\img \rho_d \subseteq \GL(F_d A)$ for a suitable $d$ as in \prref{autfl.alg.grp}, then
$$
\HH' \cap \ker \ga = \bigl\{ f \in \img \rho_d \mid (f-\id) \bigl( \bigoplus_{k=0}^n A_k \bigr) \subseteq \bigoplus_{k=0}^{n-1} A_k \,, \;\; \forall\, n \in [0,d] \bigr\}.
$$
It follows that $\HH' \cap \ker \ga$ is unipotent. Since algebraic tori have no nontrivial unipotent subgroups, $\HH' \cap \ker \ga = \{\id\}$. Therefore we conclude from $\ga(\HH') = \ga(\HH)$ that $\HH' = \HH$, as desired.
\end{proof}

We illustrate the use of the above results to recapture the twisted commutation invariant of a quantum matrix algebra from the general twist invariant of the algebra.

\bex{OlpMn3}
Let $A := \OlpMn$ as in \exref{OlpMn}, and equip $A$ with the natural $\Zset$-grading in which all the generators $X_{ij}$ have degree $1$. The first stage of the corresponding filtration on $A$ is the subspace $V := \KK + \bigoplus_{i,j=1}^n \KK X_{ij}$, which generates $A$ as a $\KK$-algebra. By Observations \ref{autfl-autgr}, $\autgr A$ becomes an algebraic group isomorphic by the restriction morphism $\rho_1$ to the closed subgroup $\AA := \rho_1(\autgr A)$ of $\GL(V)$.

Now view $A$ as a CGL extension as in \exref{OlpMn2}. If we express the maps in $\GL(V)$ as matrices with respect to the basis $1,X_{11}, X_{12}, \dots, X_{nn}$, then the action of $\HH$ on $A$ corresponds to the morphism $\Delta : \HH \rightarrow \AA$ given by the rule
\begin{multline*}
\Delta(\al_1,\dots,\al_n, \be_1,\dots, \be_n) =  \\
\text{diag}(1,\al_1\be_1, \al_1\be_2, \dots, \al_1 \be_n, \al_2\be_1, \dots, \al_n \be_1, \dots, \al_n \be_n).
\end{multline*}
Our chosen basis vectors $1,X_{11}, X_{12}, \dots, X_{nn}$ for $V$ are $\Delta(\HH)$-eigenvectors with distinct $\Delta(\HH)$-eigenvalues, and so the $\Delta(\HH)$-eigenspaces of $V$ are all $1$-dimensional. If $\HH'$ is a subtorus of $\AA$ containing $\Delta(\HH)$, all the $\Delta(\HH)$-eigenspaces of $V$ must be invariant under $\HH'$, due to the commutativity of $\HH'$. Consequently, each of  $1,X_{11}, X_{12}, \dots, X_{nn}$ must be a $\Delta(\HH)$-eigenvector. It is readily checked that the only diagonal matrices in $\AA$ are those in $\Delta(\HH)$, so $\HH' = \Delta(\HH)$. Therefore $\Delta(\HH)$ is a maximal torus of $\AA$, identified via $\rho_1^{-1}$ with a maximal torus of $\autgr A$. By \leref{maxtor.autgr}, $\rho_1^{-1}(\Delta(\HH))$ is also a maximal torus of $\autfl A$. Taking account of \eqref{ADtw.OlpMn}, we conclude that
\begin{equation}  \label{tw.OlpMn}
\tw(A) = \tw_{X(\rho_1^{-1}(\Delta(\HH)))}(A) = \tw_{\xh}(A) = \langle \la \rangle.
\end{equation}
\eex

\sectionnew{Stability of the $\AD$-invariant}
This section establishes stability properties of the $\AD$-invariant. The first result is that the invariant remains unchanged under polynomial extensions. 
The second result establishes invariance under general tensoring with commutative algebras $C$ under natural assumptions on $C$.
\subsection{Polynomial stability}
\bpr{AD[x]} If $A$ is a semiprime right Goldie $\KK$-algebra, then the polynomial ring $A[x]$ is semiprime right Goldie and
\[
\AD(A[x]) = \AD(A).
\]
\epr

\begin{proof} That $A[x]$ is semiprime right Goldie is well known (e.g., \cite[Corollary 11.19]{Lam}).  Of course, all regular elements of $A$ are regular in $A[x]$, from which it is clear that $\AD(A) \subseteq \AD(A[x])$. It remains to establish the reverse inclusion.

By Observation \ref{ADobs}(1), $\AD(A) = \AD(\Fract A)$ and $\AD(A[x]) = \AD((\Fract A)[x])$, so there is no loss of generality in assuming that $A$ is semisimple. Then, due to \eqref{ADprod}, we lose no generality in assuming that $A = M_n(D)$ for some $n \in \Zpos$ and some division $\KK$-algebra $D$.
We may identify
$$
A((x)) = M_{n}\bigl( D((x)) \bigr),
$$
where $A((x))$ and $D((x))$ denote the algebras of Laurent series in $x$ over $A$ and $D$, respectively. Since $D((x))$ is a division algebra, it follows that $A((x))$ is semisimple. Moreover, the inclusion map $D[x] \rightarrow D((x))$ extends uniquely to an injective homomorphism $\Fract D[x] \rightarrow D((x))$, from which we obtain a natural embedding
$$
\Fract A[x] \equiv M_{n}(\Fract D[x])  \longrightarrow A((x)).
$$
Consequently, it suffices to show that $\AD \bigl( A((x)) \bigr) \subseteq \AD(A)$.

Any scalar $\alpha \in \AD \bigl( A((x)) \bigr)$ has the form $\alpha = u_1v_1u_1^{-1}v_1^{-1} \cdots u_nv_nu_n^{-1}v_n^{-1}$ for some units $u_i,v_i \in A((x))$. If $a_i,b_i \in A$ denote the leading coefficients of $u_i$ and $v_i$, respectively (i.e., the coefficients on the lowest powers of $x$ appearing in these series), then $u_iv_iu_i^{-1}v_i^{-1}$ is a power series with constant term $a_ib_ia_i^{-1}b_i^{-1}$. The product of the $u_iv_iu_i^{-1}v_i^{-1}$ is thus a power series whose constant term lies in $\GL(A)'$. Therefore $\alpha \in \GL_1(A)' \cap \kx = \AD(A)$, as desired.
\end{proof}

\prref{AD[x]} immediately implies that various polynomial algebras cannot be isomorphic. In particular:

\bco{polyqaffnoniso}
Let $m,n \in \Zpos$, and let $\bfp \in M_m(\kx)$, $\bfq \in M_n(\kx)$ be multiplicatively skew-symmetric matrices. If $\langle \bfp\rangle \ne \langle \bfq \rangle$, then the polynomial algebras
$$
\OO_\bfp(\KK^m)[x_1,\dots,x_s] \not\cong \OO_\bfq(\KK^n)[x_1,\dots,x_t]
$$
for all $s,t\in \Znn$. 

More generally, for all pairs of quantum cluster algebras $\AA_1$ and $\AA_2$ {\em{(}}and pairs of upper quantum cluster algebras 
$\UU_1$ and $\UU_2${\em{)}} we have 
\[
\AA_1[x_1,\dots,x_s] \not\cong \AA_2[x_1,\dots,x_t] \quad \mbox{if} \quad q(\AA_1) \neq q(\AA_2)
\]
and
\[
\UU_1[x_1,\dots,x_s] \not\cong \UU_2[x_1,\dots,x_t] \quad \mbox{if} \quad q(\UU_1) \neq q(\UU_2)
\]
for all $s,t\in \Znn$. 
\eco
\subsection{General stability under tensoring with commutative algebras}
\ble{nonhomfract}
Let $R$ be a ring and $X\subseteq R$ a right denominator set. Suppose $a_0,\dots,a_n\in R$ and $x_0,\dots,x_n \in X$ such that
\begin{equation}  \label{fracprod}
a_0x_0^{-1} = (a_1x_1^{-1}) (a_2x_2^{-1}) \cdots (a_nx_n^{-1})
\end{equation}
in $RX^{-1}$. 

{\rm(a)} There exist $b_1,\dots,b_n\in R$ and $y_1=1, y_2,\dots,y_n,z \in X$ such that
\begin{equation}  \label{fracprodcond}
\begin{aligned}
a_1 b_2 \cdots b_n z &= a_0 b_1  &&\qquad x_n y_n z = x_0 b_1  \\
x_{i-1} y_{i-1} b_i &= a_i y_i &&\qquad \forall\, i \in [2,n].
\end{aligned}
\end{equation}

{\rm(b)} Suppose $\phi : R \rightarrow S$ is a ring homomorphism which maps $x_0,\dots, x_n, y_2,\dots, y_n, z$ to units of $S$. Then
\begin{equation}  \label{phifrac}
\phi(a_0) \phi(x_0)^{-1} = \phi(a_1) \phi(x_1)^{-1} \phi(a_2) \phi(x_2)^{-1} \cdots \phi(a_n) \phi(x_n)^{-1}.
\end{equation}
\ele

\begin{proof} (a) By induction on $n$, there exist $b_2,\dots,b_n\in R$ and $y_2,\dots,y_n \in X$ such that $x_{i-1} b_i = a_i y_i$ for $i\in [2,n]$ and
$$
(a_1x_1^{-1}) (a_2x_2^{-1}) \cdots (a_nx_n^{-1}) = (a_1 b_2 \cdots b_n) (x_n y_n)^{-1}.
$$
There also exist $b_1\in R$ and $z\in X$ such that $x_n y_n z = x_0 b_1$. Incorporating \eqref{fracprod} and multiplying on the right by $x_0 b_1$ yields $a_0 b_1 = a_1 b_2 \cdots b_n z$.

(b) Apply $\phi$ to \eqref{fracprodcond} and unwind.
\end{proof}

In the following, we let $\udim S$ denote the right uniform dimension of a ring $S$, i.e., the uniform dimension of the module $S_S$.

\bpr{patchapplic}
Let $R$ be a semiprime right noetherian ring and $\PP$ a set of prime ideals of $R$ such that $\bigcap \PP = 0$ and there is a finite upper bound on $\{ \udim R/P \mid P \in \PP \}$. For any regular element $c\in R$, there is some $P \in \PP$ such that $c$ is regular modulo $\PP$.
\epr

\begin{proof} There are only finitely many distinct minimal prime ideals in $R$ (e.g., \cite[Theorem 3.4]{GW}), say $P_1,\dots,P_n$, and $c$ is regular modulo each $P_i$ (e.g., \cite[Lemma 7.4]{GW}). If $\PP_i := \{ P \in \PP \mid P \supseteq P_i \}$, then $\bigcap_{i=1}^n \bigcap \PP_i = \bigcap \PP = 0$. Consequently, $\bigcap \PP_j \subseteq P_1$ for some $j$. Since $P_j \subseteq \bigcap \PP_j$, we must have $j=1$ and $\bigcap \PP_1 = P_1$.

Since $\{ \udim R/P \mid P \in \PP_1 \}$ is bounded above, \cite[Proposition 16.10]{GW} implies that there is a neighborhood $\UU$ of $P_1$ in the patch topology on $\Spec R$ such that $c$ is regular modulo $P$ for all $P \in \PP_1\cap \UU$. The discussion in \cite[pp.~275-6]{GW} shows that $\UU$ contains a patch neighborhood of the form 
$$\UU' = \{ Q \in \Spec R \mid Q \supseteq P_1 \; \text{but} \; Q \nsupseteq I \},$$
for some ideal $I$ properly containing $P_1$. Since $\bigcap \PP_1 = P_1$, we conclude that $\PP_1 \cap \UU'$ is nonempty. Any $P \in \PP_1 \cap \UU'$ satisfies the requirements of the proposition.
\end{proof}

\bth{ADtensorC}
Assume $\KK = \kbar$, and let $C$ be a reduced {\rm(}i.e., semiprime\/{\rm)} commutative $\KK$-algebra. If $A$ is a semiprime right Goldie $\KK$-algebra and $A\otimes_\KK C$ is semiprime right Goldie, then
\begin{equation}  \label{ADtensor}
\AD(A\otimes_\KK C) = \AD(A).
\end{equation}
\eth

\begin{proof} We just need to show that $\AD(A\otimes_\KK C) \subseteq \AD(A)$, since the reverse inclusion is clear.

If $X$ is the set of regular elements of $A$, then $X$ is a right denominator set in $A$ and $X\otimes 1$ is a right denominator set of regular elements in $A\otimes_\KK C$. Consequently, $\AD(A) = \AD(AX^{-1})$ and $\AD(A\otimes_\KK C) = \AD((A\otimes_\KK C)(X\otimes1)^{-1})$, so there is no loss of generality in replacing $A$ by $AX^{-1}$. Thus, we may now assume that $A$ is semisimple. In view of \eqref{ADprod}, we may further assume that $A$ is simple artinian.

Now $C$ is a directed union of its affine subalgebras $C'$, so $A\otimes_\KK C$ is a directed union of the subalgebras $A\otimes_\KK C'$. Since $A\otimes_\KK C$ is a central extension of $A\otimes_\KK C'$, semiprimeness descends from $A\otimes_\KK C$ to $A\otimes_\KK C'$. Moreover, $A\otimes_\KK C'$ is a homomorphic image of a polynomial ring $A[x_1,\dots,x_n]$, and so $A\otimes_\KK C'$ is right noetherian, hence right Goldie. Note that any element of $\AD(A\otimes_\KK C)$ lies in $\AD(A\otimes_\KK C')$ for a suitable $C'$. Thus, it suffices to prove that $\AD(A\otimes_\KK C') \subseteq \AD(A)$ for all $C'$. Consequently, without loss of generality, we may assume that $C$ is affine and $A\otimes_\KK C$ is right noetherian.

Let $\Phi$ denote the set of characters $\phi : C \rightarrow \KK$ and $\Psi$ the set of $\KK$-algebra homomorphisms of the form
$$
\id_A \otimes \phi : A\otimes_\KK C \longrightarrow A,
$$
for $\phi \in \Phi$. Since $C$ is affine and reduced and $\KK$ is algebraically closed, $\bigcap_{\phi \in \Phi} \ker \phi = 0$, from which it follows that $\bigcap_{\phi\in \Phi} \ker(\id_A \otimes \phi) = 0$. In other words, the collection
$$
\PP := \{ \ker \psi \mid \psi \in \Psi \}
$$
of maximal ideals of $A\otimes_\KK C$ has zero intersection. Note that the algebras $(A\otimes_\KK C)/P$, for $P \in \PP$, all have the same uniform dimension, namely $\udim A$.

Let $X$ denote the set of regular elements of $A\otimes_\KK C$. Any scalar $\alpha \in \AD(A\otimes_\KK C)$ has the form
$$
\alpha = s_1 t_1^{-1} u_1 v_1^{-1} t_1 s_1^{-1} v_1 u_1^{-1} \cdots s_n t_n^{-1} u_n v_n^{-1} t_n s_n^{-1} v_n u_n^{-1}
$$
for some $s_i,t_i,u_i,v_i \in X$. Relabel the sequences
$$
\alpha, s_1, u_1, t_1, v_1, \dots, s_n, u_n, t_n, v_n \qquad\text{and}\qquad 1, t_1, v_1, s_1, u_1, \dots, t_n, v_n, s_n, u_n,
$$
as $a_0,\dots,a_{4n}$ and $x_0,\dots,x_{4n}$, respectively. By \leref{nonhomfract}(a), there exist $b_1,\dots, b_{4n} \in R$ and $y_1=1, y_2, \dots, y_{4n}, z \in X$ such that
\begin{align*}
a_1 b_2 \cdots b_{4n} z &= a_0 b_1  &\qquad&x_{4n} y_{4n} z = x_0 b_1  \\
x_{i-1} y_{i-1} b_i &= a_i y_i  &&\forall\, i \in [2,4n].
\end{align*}
Now let $c$ be the product (in some order) of the $x_i$, $y_i$, and $z$, a regular element of $A\otimes_\KK C$. By \prref{patchapplic}, there exists $P \in \PP$ such that $c$ is regular modulo $P$. Then $P = \ker \psi$ for some $\psi  \in \Psi$, and $\psi(c)$ is a unit in the simple artinian ring $A$. Consequently, the $\psi(x_i)$, $\psi(y_i)$, and $\psi(z)$ are units in $A$. Applying \leref{nonhomfract}(b) and rewriting the result in terms of the original $s_i$, $t_i$, $u_i$, $v_i$, we obtain that the $\psi(s_i)$, $\psi(t_i)$, $\psi(u_i)$, $\psi(v_i)$ are units in $A$ and
\begin{multline*}
\alpha = \psi(s_1) \psi(t_1)^{-1} \psi(u_1) \psi(v_1)^{-1} \psi(t_1) \psi(s_1)^{-1} \psi(v_1) \psi(u_1)^{-1} \cdots  \\
\psi(s_n) \psi(t_n)^{-1} \psi(u_n) \psi(v_n)^{-1} \psi(t_n) \psi(s_n)^{-1} \psi(v_n) \psi(u_n)^{-1}.
\end{multline*}
Thus $\alpha \in \AD(A)$, as required.
\end{proof} 

\bco{CGLtensorC}
Assume $\KK = \kbar$, and let $C$ be a reduced commutative noetherian $\KK$-algebra. If $A$ is an iterated Ore extension of $\KK$, then
$$
\AD(A\otimes_\KK C) = \AD(A).
$$
\eco

\begin{proof} If $X$ is the set of regular elements of $C$, then $CX^{-1}$ is a finite direct product of field extensions $K_i$ of $\KK$. Now $1\otimes X$ is a central multiplicative set in $A\otimes_\KK C$, and $(A\otimes_\KK C)(1\otimes X)^{-1}$ is isomorphic to the finite direct product of the $\KK$-algebras $A\otimes_\KK K_i$. Each of the latter tensor products is an iterated Ore extension of $K_i$, hence a noetherian domain. Thus, $(A\otimes_\KK C)(1\otimes X)^{-1}$ is a semiprime noetherian ring, from which we conclude that $A\otimes_\KK C$ is semiprime Goldie. The result now follows from \thref{ADtensorC}.
\end{proof}

\sectionnew{Stability of the twist invariants}
In this section we prove that both twist invariants remain unchanged under polynomial extensions. 
We apply the fact for the general twist invariant to establish non-isomorphisms between twisted polynomial extensions of quantized Weyl algebras.
\subsection{Maximal tori of polynomial extensions}
First we show that the maximal tori of the automorphism groups of filtered locally finite algebras behave well under polynomial extensions. For any such algebra $A$, we give the polynomial ring $A[x]$ the natural induced filtration with
\begin{equation}  \label{polyfilt}
F_k(A[x]) := \bigoplus_{j=0}^k\, (F_j A) x^{k-j}, \;\; \forall\, k \in \Znn \,.
\end{equation}

\bpr{A[x]-maxtor} 
Let $A$ be an affine filtered locally finite algebra over an infinite field $\KK$, with $F_0 A = \KK$. For every maximal torus $\HH$ of $\autfl(A)$, the group $\HH \times \kx$ 
is a maximal torus of $\autfl(A[x])$, where $\HH$ acts trivially on $x$ and $\kx$ acts by rescaling $x$.
\epr

\begin{proof} Set $\HH^+ := \HH \times \kx$, and observe that the $\HH^+$-eigenspaces of $A[x]$ are the spaces $Ex^n$ where $E$ is an $\HH$-eigenspace of $A$ and $n \in \Znn$. In particular, $A[x]$ is a direct sum of $\HH^+$-eigenspaces.

Suppose $\HH'$ is a subtorus of $\autfl(A[x])$ containing $\HH^+$. Since $\HH'$ is abelian, it stabilizes the $\HH^+$-eigenspaces of $A[x]$, and so it stabilizes $A$. Thus, restriction to $A$ provides a group homomorphism $\rho_1 : \HH' \rightarrow \autfl(A)$, and we observe that $\rho_1$ is a morphism of algebraic groups. Now $\HH'$ also stabilizes $F_1(A[x]) = (F_1A)+ \KK x$ as well as the $\HH^+$-eigenspace $A^\HH x$, where $A^\HH$ is the fixed subalgebra of $A$ under $\HH$, so $\HH'$ stabilizes $\KK x$. Consequently, restriction to $\KK x$ provides a second algebraic group morphism $\rho_2 : \HH' \rightarrow \GL(\KK x) = \kx$. The combined morphism $\rho_1 \times \rho_2 : \HH' \rightarrow \autfl(A) \times \kx$ is injective, because the automorphisms in $\HH'$ are determined by their actions on $A \cup \KK x$. Obviously $(\rho_1\times \rho_2)(\HH^+) = \HH^+$. Hence, once we show that $(\rho_1 \times \rho_2)(\HH') = \HH^+$, we can conclude that $\HH' = \HH^+$, establishing the maximality of $\HH^+$.

Choose $n>0$ large enough so that $F_n A$ generates $A$, whence $F_n(A[x])$ generates $A[x]$. The algebraic group structure on $\autfl(A[x])$ thus comes from the isomorphism of this group with a closed subgroup of $\GL(F_n(A[x]))$. Since the group $\autfl(A[x])$ acts rationally on $A[x]$, so does $\HH'$. This means $\HH'$ acts diagonalizably on $F_n(A[x])$. Now $F_n(A[x])$ is $\HH^+$-stable, and some of its $\HH^+$-eigenspaces sum to $F_nA$, from which we see that $\HH'$ stabilizes $F_nA$. Consequently, $\HH'$ acts diagonalizably on $F_nA$. It follows that $\rho_1(\HH')$ is a diagonalizable algebraic group. On the other hand, it is connected, being a continuous image of a connected group. Thus, $\rho_1(\HH')$ is a torus. Since $\rho_1(\HH') \supseteq \HH$, we obtain $\rho_1(\HH') = \HH$ from the maximality of $\HH$. Therefore $(\rho_1 \times \rho_2)(\HH') \subseteq \HH^+$ and consequently $(\rho_1 \times \rho_2)(\HH') = \HH^+$, as desired. 
 \end{proof}
 
\subsection{Stability of the twist invariants}

\bth{twA[x]} {\em{(a)}} Let $A$ be a $\Gamma$-graded semiprime right Goldie $\KK$-algebra for an abelian group $\Gamma$. 
Then the polynomial algebra $A[x]$ is also semiprime right Goldie and
\[
\tw_{\Gamma \times \Zset}(A[x]) = \tw_\Gamma(A),
\]
where the grading of $A[x]$ is obtained by extending that of $A$ by placing $x$ in degree $(0,1)$. 

{\em{(b)}} Let $A$ be an affine filtered locally finite algebra over an infinite field $\KK$, with $F_0 A = \KK$.
Assume in addition that $A$ is a semiprime right Goldie algebra. Then the polynomial algebra $A[x]$ has the same properties 
with respect to the filtration \eqref{polyfilt} and 
\[
\tw(A[x]) = \tw(A).
\]
\eth

\begin{proof} (a) The fact that $A[x]$ is also semiprime right Goldie is discussed in \prref{AD[x]}. 
Each $\cb \in Z^2(\Ga, \kx)$ gives rise to 
\[
\cb^+ \in Z^2( \Ga \times \Zset, \kx) \quad \mbox{defined by} \quad 
\cb^+( (a_1, n_1), (a_2, n_2)):= \cb(a_1, a_2)
\]
for $a_i \in \Ga$, $n_i \in \Zset$. Obviously, 
\[
(A[x])^{\cb^+} \cong A^\cb[x],
\]
and by \prref{AD[x]}, $\AD(A^\cb[x]) = \AD(A^\cb)$. Therefore, 
\begin{align*}
\tw_\Ga(A) &= \bigcap_{\cb \in Z^2(\Ga,\kx)} \AD(A^\cb) = \bigcap_{\cb \in Z^2(\Ga,\kx)} \AD((A[x])^{\cb^+}) \supseteq
\\
&\supseteq \bigcap_{\db \in Z^2(\Ga \times \Zset,\kx)} \AD((A[x])^\db) = \tw_{\Ga \times \Zset} (A[x]).
\end{align*}
In the opposite direction, for $\db \in Z^2(\Ga \times \Zset,\kx)$, we have the restricted cocycle 
\[
\ol{\db} := \db|_{\Ga \times \Ga} \in Z^2(\Ga,\kx).
\]
Since $\Fract(A^{\ol{\db}}) \subset \Fract((A[x])^\db)$, 
\[
\AD((A[x])^\db) \supseteq \AD ( A^{\ol{\db}}).
\]
Thus,
\begin{align*}
\tw_{\Ga \times \Zset}(A[x]) &= \bigcap_{\db \in Z^2(\Ga \times \Zset,\kx)} \AD((A[x])^\db) 
\supseteq \bigcap_{\db \in Z^2(\Ga \times \Zset,\kx)} \AD(A^{\ol{\db}}) \supseteq 
\\
&\supseteq  \bigcap_{\cb \in Z^2(\Ga,\kx)} \AD(A^\cb)  = \tw_\Ga (A).
\end{align*}
Combining the two inclusions above gives $\tw_{\Ga \times \Zset}(A[x]) = \tw_\Ga (A)$.

(b) The properties of $A[x]$ are discussed in the proofs of Propositions \ref{pAD[x]} and \ref{pA[x]-maxtor}. 
Let $\HH$ be a maximal torus of $\autfl (A) $. By \prref{A[x]-maxtor}, $\HH^+ := \HH \times \kx$ is a maximal torus of $\autfl( A[x])$, 
where $\HH$ acts trivially on $x$ and $\kx$ acts by rescaling $x$. We identify $X(\HH^+) \cong \xh \times \Zset$ in the canonical way. 
Applying the first part of the theorem and taking into account \leref{tw.indepH}, we obtain
\[
\tw (A[x]) = \tw_{X(\HH^+)}(A[x]) = \tw_{\xh \times \Zset}(A[x]) = 
\tw_{\xh}(A) = \tw(A),
\]
which completes the proof of the theorem. 
\end{proof}

\bex{OlpMn4}
Let $A := \OlpMn$ as in Examples \ref{eOlpMn} and \ref{eOlpMn3}, treated as a $\Zset$-graded algebra with $\deg X_{ij} = 1$ for all $i$, $j$. Extend this grading to any polynomial algebra $A[x_1,\dots,x_s]$ so that $\deg x_j = 1$ for $j \in [1,s]$. In view of \thref{twA[x]}(b) and equation \eqref{tw.OlpMn}, we obtain
\begin{equation}  \label{tw.OlpMn[x]}
\tw(A[x_1,\dots,x_s]) = \langle \la \rangle, \;\; \forall\, \text{polynomial algebras} \;\; A[x_1,\dots,x_s].
\end{equation}
\eex

\bex{qWeyl3}
Let $A := A_n^{Q,P}(\KK)$ as in \exref{qWeyl}, treated as a $\Zset$-graded algebra with $\deg x_i = 1$ and $\deg y_i = -1$ for all $i \in [1,n]$. This grading is not  locally finite, but $A$ has a  locally finite filtration based on the generating subspace $V: = \KK \oplus \bigoplus_{i=1}^n (\KK x_i \oplus \KK y_i)$, where we take $F_kA = V^k$ for $k \in \Znn$. By Observations \ref{autfl-autgr}, $\autfl A$ becomes an algebraic group isomorphic by the restriction morphism $\rho_1$ to the closed subgroup $\AA := \img \rho_1$ of $\GL(V)$.

Let the torus $\HH := (\kx)^n$ act on $A$ as in \eqref{Hact.qWeyl} (but here we allow the $q_i$ to possibly be roots of unity), and note that the corresponding $\xh$-grading on $A$ coincides with the $\Zset^n$-grading of \exref{qWeyl}. Obviously $\HH$ acts by filtered automorphisms, and we may identify $\HH$ with a closed subgroup of $\autfl A$. Note that the $\HH$-eigenspaces of $V$ are all $1$-dimensional. Hence, it follows as in \exref{OlpMn3} that $\HH$ is a maximal torus of $\autfl A$. By \eqref{twGa.qWeyl}, we have
\begin{equation}  \label{tw.qWeyl}
\tw(A) = \tw_{\xh}(A) = \langle q_1, \dots, q_n \rangle,
\end{equation}
and \thref{twA[x]}(b) thus implies
$$
\tw(A[X_1,\dots, X_s]) = \langle q_1, \dots, q_n \rangle,  \;\; \forall\, \text{polynomial algebras} \;\; A[X_1,\dots,X_s],
$$
where $A[X_1,\dots,X_s]$ is filtered by iterating Eqn.~\eqref{polyfilt}:
$$
F_k (A[X_1,\dots,X_s]) := \bigoplus_{j_0+ j_1+ \cdots+ j_s \le k} (F_{j_0} A) X_1^{j_1} \cdots X_s^{j_s}.
$$
Taking all polynomial algebras over quantized Weyl algebras to be filtered in this manner, we conclude that
\begin{enumerate}
\item[] Given $A_n^{Q,P}(\KK)$ and $A_{n'}^{Q',P'}(\KK)$ with $\langle q_1, \dots, q_n \rangle \ne \langle q'_1, \dots, q'_n \rangle$, no cocycle twisted 
polynomial algebra 
\[
\big( A_n^{Q,P}(\KK)[X_1, \dots, X_s] \big)^{\cb}
\]
is isomorphic, as a filtered $\KK$-algebra, to any cocycle twisted polynomial algebra 
\[
\big( A_{n'}^{Q',P'}(\KK)[X_1, \dots, X_t] \big)^{\cb'}
\]
for any two 2-cocycles $\cb \in Z^2(\Zset^{m},\kx)$ and $\cb' \in Z^2(\Zset^{m'},\kx)$. The twists 
are performed with respect to any gradings of the two polynomial algebras that respect the filtrations
of these algebras.
\end{enumerate}
In the last step we made use of the invariance property of the general twist invariant from \thref{invar-tw}. 
\eex

We finish by noting that recently Bell and Zhang \cite{BeZh} investigated Zariski cancellative algebras, which are $\KK$-algebras $A$ such that  
$A[x] \cong B[x]$ implies $A \cong B$ for any $\KK$-algebra $B$. They proved that this is the case if $A$ is ${\mathrm{LND}}$-rigid or if the discriminant of $A$ over 
its center exists and is cancellative. Currently, there are not many families of algebras that are known to have these properties and 
even quantum matrices of size $n > 3$ are conjectured not to have the properties.
Our stable twist invariant can be explicitly computed from Theorems \ref{ttwGa.sandwich} and \ref{ttwGa.cgl} and used to detect non-isomorphisms 
\[
A[x] \ncong B[x]
\]
of filtered algebras when discriminant methods are not applicable.

\end{document}